
\ifx\shlhetal\undefinedcontrolsequence\let\shlhetal\relax\fi
\def\fmtname{AmS-TeX}

\def\fmtversion{2.1}
\catcode`\@=11
\ifx\amstexloaded@\relax\catcode`\@=\active
  \endinput\else\let\amstexloaded@\relax\fi
\newlinechar=`\^^J
\def\W@{\immediate\write\sixt@@n}
\def\CR@{\W@{^^J\fmtname - Version \fmtversion^^J%
COPYRIGHT 1985, 1990, 1991 - AMERICAN MATHEMATICAL SOCIETY^^J%
Use of this macro package is not restricted provided^^J%
each use is acknowledged upon publication.^^J}}
\CR@ \everyjob{\CR@}
\message{Loading definitions for}
\message{misc utility macros,}
\toksdef\toks@@=2
\long\def\rightappend@#1\to#2{\toks@{\\{#1}}\toks@@
 =\expandafter{#2}\xdef#2{\the\toks@@\the\toks@}\toks@{}\toks@@{}}
\def\alloclist@{}
\newif\ifalloc@
\def\showallocations{{\def\\{\immediate\write\m@ne}\alloclist@}\alloc@true}
\def\alloc@#1#2#3#4#5{\global\advance\count1#1by\@ne
 \ch@ck#1#4#2\allocationnumber=\count1#1
 \global#3#5=\allocationnumber
 \edef\next@{\string#5=\string#2\the\allocationnumber}%
 \expandafter\rightappend@\next@\to\alloclist@}
\newcount\count@@
\newcount\count@@@
\def\FN@{\futurelet\next}
\def\DN@{\def\next@}
\def\DNii@{\def\nextii@}
\def\RIfM@{\relax\ifmmode}
\def\RIfMIfI@{\relax\ifmmode\ifinner}
\def\setboxz@h{\setbox\z@\hbox}
\def\wdz@{\wd\z@}
\def\boxz@{\box\z@}
\def\setbox@ne{\setbox\@ne}
\def\wd@ne{\wd\@ne}
\def\iterate{\body\expandafter\iterate\else\fi}
\def\err@#1{\errmessage{AmS-TeX error: #1}}
\newhelp\defaulthelp@{Sorry, I already gave what help I could...^^J
Maybe you should try asking a human?^^J
An error might have occurred before I noticed any problems.^^J
``If all else fails, read the instructions.''}
\def\Err@{\errhelp\defaulthelp@\err@}
\def\eat@#1{}
\def\in@#1#2{\def\in@@##1#1##2##3\in@@{\ifx\in@##2\in@false\else\in@true\fi}%
 \in@@#2#1\in@\in@@}
\newif\ifin@
\def\space@.{\futurelet\space@\relax}
\space@. %
\newhelp\athelp@
{Only certain combinations beginning with @ make sense to me.^^J
Perhaps you wanted \string\@\space for a printed @?^^J
I've ignored the character or group after @.}
{\catcode`\~=\active 
 \lccode`\~=`\@ \lowercase{\gdef~{\FN@\at@}}}
\def\at@{\let\next@\at@@
 \ifcat\noexpand\next a\else\ifcat\noexpand\next0\else
 \ifcat\noexpand\next\relax\else
   \let\next\at@@@\fi\fi\fi
 \next@}
\def\at@@#1{\expandafter
 \ifx\csname\space @\string#1\endcsname\relax
  \expandafter\at@@@ \else
  \csname\space @\string#1\expandafter\endcsname\fi}
\def\at@@@#1{\errhelp\athelp@ \err@{\Invalid@@ @}}
\def\atdef@#1{\expandafter\def\csname\space @\string#1\endcsname}
\newhelp\defahelp@{If you typed \string\define\space cs instead of
\string\define\string\cs\space^^J
I've substituted an inaccessible control sequence so that your^^J
definition will be completed without mixing me up too badly.^^J
If you typed \string\define{\string\cs} the inaccessible control sequence^^J
was defined to be \string\cs, and the rest of your^^J
definition appears as input.}
\newhelp\defbhelp@{I've ignored your definition, because it might^^J
conflict with other uses that are important to me.}
\def\define{\FN@\define@}
\def\define@{\ifcat\noexpand\next\relax
 \expandafter\define@@\else\errhelp\defahelp@                               
 \err@{\string\define\space must be followed by a control
 sequence}\expandafter\def\expandafter\nextii@\fi}                          
\def\undefined@@@@@@@@@@{}
\def\preloaded@@@@@@@@@@{}
\def\next@@@@@@@@@@{}
\def\define@@#1{\ifx#1\relax\errhelp\defbhelp@                              
 \err@{\string#1\space is already defined}\DN@{\DNii@}\else
 \expandafter\ifx\csname\expandafter\eat@\string                            
 #1@@@@@@@@@@\endcsname\undefined@@@@@@@@@@\errhelp\defbhelp@
 \err@{\string#1\space can't be defined}\DN@{\DNii@}\else
 \expandafter\ifx\csname\expandafter\eat@\string#1\endcsname\relax          
 \global\let#1\undefined\DN@{\def#1}\else\errhelp\defbhelp@
 \err@{\string#1\space is already defined}\DN@{\DNii@}\fi
 \fi\fi\next@}

\def\predefine#1#2{\let#1#2}
\def\undefine#1{\let#1\undefined}
\message{page layout,}
\newdimen\captionwidth@
\captionwidth@\hsize
\advance\captionwidth@-1.5in
\def\pagewidth#1{\hsize#1\relax
 \captionwidth@\hsize\advance\captionwidth@-1.5in}
\def\pageheight#1{\vsize#1\relax}
\def\hcorrection#1{\advance\hoffset#1\relax}
\def\vcorrection#1{\advance\voffset#1\relax}
\message{accents/punctuation,}

\let\graveaccent\`
\let\acuteaccent\'
\let\tildeaccent\~
\let\hataccent\^
\let\underscore\_
\let\B\=
\let\D\.
\let\ic@\/
\def\/{\unskip\ic@}
\def\textfonti{\the\textfont\@ne}
\def\t#1#2{{\edef\next@{\the\font}\textfonti\accent"7F \next@#1#2}}
\def~{\unskip\nobreak\ \ignorespaces}
\def\.{.\spacefactor\@m}
\atdef@;{\leavevmode\null;}
\atdef@:{\leavevmode\null:}
\atdef@?{\leavevmode\null?}
\edef\@{\string @}
\def\relaxnext@{\let\next\relax}
\atdef@-{\relaxnext@\leavevmode
 \DN@{\ifx\next-\DN@-{\FN@\nextii@}\else
  \DN@{\leavevmode\hbox{-}}\fi\next@}%
 \DNii@{\ifx\next-\DN@-{\leavevmode\hbox{---}}\else
  \DN@{\leavevmode\hbox{--}}\fi\next@}%
 \FN@\next@}
\def\srdr@{\kern.16667em}
\def\drsr@{\kern.02778em}
\def\sldl@{\drsr@}
\def\dlsl@{\srdr@}
\atdef@"{\unskip\relaxnext@
 \DN@{\ifx\next\space@\DN@. {\FN@\nextii@}\else
  \DN@.{\FN@\nextii@}\fi\next@.}%
 \DNii@{\ifx\next`\DN@`{\FN@\nextiii@}\else
  \ifx\next\lq\DN@\lq{\FN@\nextiii@}\else
  \DN@####1{\FN@\nextiv@}\fi\fi\next@}%
 \def\nextiii@{\ifx\next`\DN@`{\sldl@``}\else\ifx\next\lq
  \DN@\lq{\sldl@``}\else\DN@{\dlsl@`}\fi\fi\next@}%
 \def\nextiv@{\ifx\next'\DN@'{\srdr@''}\else
  \ifx\next\rq\DN@\rq{\srdr@''}\else\DN@{\drsr@'}\fi\fi\next@}%
 \FN@\next@}

\def\textfontii{\the\textfont\tw@}
\def\lbrace@{\delimiter"4266308 }
\def\rbrace@{\delimiter"5267309 }
\def\{{\RIfM@\lbrace@\else{\textfontii f}\spacefactor\@m\fi}
\def\}{\RIfM@\rbrace@\else
 \let\@sf\empty\ifhmode\edef\@sf{\spacefactor\the\spacefactor}\fi
 {\textfontii g}\@sf\relax\fi}
\let\lbrace\{
\let\rbrace\}
\def\AmSTeX{{\textfontii A\kern-.1667em%
  \lower.5ex\hbox{M}\kern-.125emS}-\TeX}
\message{line and page breaks,}
\def\vmodeerr@#1{\Err@{\string#1\space not allowed between paragraphs}}
\def\mathmodeerr@#1{\Err@{\string#1\space not allowed in math mode}}
\def\linebreak{\RIfM@\mathmodeerr@\linebreak\else
 \ifhmode\unskip\unkern\break\else\vmodeerr@\linebreak\fi\fi}

\newskip\saveskip@
\def\allowlinebreak{\RIfM@\mathmodeerr@\allowlinebreak\else
 \ifhmode\saveskip@\lastskip\unskip
 \allowbreak\ifdim\saveskip@>\z@\hskip\saveskip@\fi
 \else\vmodeerr@\allowlinebreak\fi\fi}
\def\nolinebreak{\RIfM@\mathmodeerr@\nolinebreak\else
 \ifhmode\saveskip@\lastskip\unskip
 \nobreak\ifdim\saveskip@>\z@\hskip\saveskip@\fi
 \else\vmodeerr@\nolinebreak\fi\fi}
\def\newline{\relaxnext@
 \DN@{\RIfM@\expandafter\mathmodeerr@\expandafter\newline\else
  \ifhmode\ifx\next\par\else
  \expandafter\unskip\expandafter\null\expandafter\hfill\expandafter\break\fi
  \else
  \expandafter\vmodeerr@\expandafter\newline\fi\fi}%
 \FN@\next@}
\def\dmatherr@#1{\Err@{\string#1\space not allowed in display math mode}}
\def\nondmatherr@#1{\Err@{\string#1\space not allowed in non-display math
 mode}}
\def\onlydmatherr@#1{\Err@{\string#1\space allowed only in display math mode}}
\def\nonmatherr@#1{\Err@{\string#1\space allowed only in math mode}}
\def\mathbreak{\RIfMIfI@\break\else
 \dmatherr@\mathbreak\fi\else\nonmatherr@\mathbreak\fi}
\def\nomathbreak{\RIfMIfI@\nobreak\else
 \dmatherr@\nomathbreak\fi\else\nonmatherr@\nomathbreak\fi}
\def\allowmathbreak{\RIfMIfI@\allowbreak\else
 \dmatherr@\allowmathbreak\fi\else\nonmatherr@\allowmathbreak\fi}
\def\pagebreak{\RIfM@
 \ifinner\nondmatherr@\pagebreak\else\postdisplaypenalty-\@M\fi
 \else\ifvmode\removelastskip\break\else\vadjust{\break}\fi\fi}
\def\nopagebreak{\RIfM@
 \ifinner\nondmatherr@\nopagebreak\else\postdisplaypenalty\@M\fi
 \else\ifvmode\nobreak\else\vadjust{\nobreak}\fi\fi}
\def\nonvmodeerr@#1{\Err@{\string#1\space not allowed within a paragraph
 or in math}}
\def\vnonvmode@#1#2{\relaxnext@\DNii@{\ifx\next\par\DN@{#1}\else
 \DN@{#2}\fi\next@}%
 \ifvmode\DN@{#1}\else
 \DN@{\FN@\nextii@}\fi\next@}
\def\newpage{\vnonvmode@{\vfill\break}{\nonvmodeerr@\newpage}}
\def\smallpagebreak{\vnonvmode@\smallbreak{\nonvmodeerr@\smallpagebreak}}
\def\medpagebreak{\vnonvmode@\medbreak{\nonvmodeerr@\medpagebreak}}
\def\bigpagebreak{\vnonvmode@\bigbreak{\nonvmodeerr@\bigpagebreak}}
\def\NoBlackBoxes{\global\overfullrule\z@}
\def\BlackBoxes{\global\overfullrule5\p@}
\def\Invalid@#1{\def#1{\Err@{\Invalid@@\string#1}}}
\def\Invalid@@{Invalid use of }
\message{figures,}
\Invalid@\caption
\Invalid@\captionwidth
\newdimen\smallcaptionwidth@
\def\topspace{\mid@false\ins@}
\def\midspace{\mid@true\ins@}
\newif\ifmid@
\def\captionfont@{}
\def\ins@#1{\relaxnext@\allowbreak
 \smallcaptionwidth@\captionwidth@\gdef\thespace@{#1}%
 \DN@{\ifx\next\space@\DN@. {\FN@\nextii@}\else
  \DN@.{\FN@\nextii@}\fi\next@.}%
 \DNii@{\ifx\next\caption\DN@\caption{\FN@\nextiii@}%
  \else\let\next@\nextiv@\fi\next@}%
 \def\nextiv@{\vnonvmode@
  {\ifmid@\expandafter\midinsert\else\expandafter\topinsert\fi
   \vbox to\thespace@{}\endinsert}
  {\ifmid@\nonvmodeerr@\midspace\else\nonvmodeerr@\topspace\fi}}%
 \def\nextiii@{\ifx\next\captionwidth\expandafter\nextv@
  \else\expandafter\nextvi@\fi}%
 \def\nextv@\captionwidth##1##2{\smallcaptionwidth@##1\relax\nextvi@{##2}}%
 \def\nextvi@##1{\def\thecaption@{\captionfont@##1}%
  \DN@{\ifx\next\space@\DN@. {\FN@\nextvii@}\else
   \DN@.{\FN@\nextvii@}\fi\next@.}%
  \FN@\next@}%
 \def\nextvii@{\vnonvmode@
  {\ifmid@\expandafter\midinsert\else
  \expandafter\topinsert\fi\vbox to\thespace@{}\nobreak\smallskip
  \setboxz@h{\noindent\ignorespaces\thecaption@\unskip}%
  \ifdim\wdz@>\smallcaptionwidth@\centerline{\vbox{\hsize\smallcaptionwidth@
   \noindent\ignorespaces\thecaption@\unskip}}%
  \else\centerline{\boxz@}\fi\endinsert}
  {\ifmid@\nonvmodeerr@\midspace
  \else\nonvmodeerr@\topspace\fi}}%
 \FN@\next@}
\message{comments,}
\def\newcodes@{\catcode`\\12\catcode`\{12\catcode`\}12\catcode`\#12%
 \catcode`\%12\relax}
\def\oldcodes@{\catcode`\\0\catcode`\{1\catcode`\}2\catcode`\#6%
 \catcode`\%14\relax}
\def\comment{\newcodes@\endlinechar=10 \comment@}
{\lccode`\0=`\\
\lowercase{\gdef\comment@#1^^J{\comment@@#10endcomment\comment@@@}%
\gdef\comment@@#10endcomment{\FN@\comment@@@}%
\gdef\comment@@@#1\comment@@@{\ifx\next\comment@@@\let\next\comment@
 \else\def\next{\oldcodes@\endlinechar=`\^^M\relax}%
 \fi\next}}}
\def\pr@m@s{\ifx'\next\DN@##1{\prim@s}\else\let\next@\egroup\fi\next@}
\def\prime{{\null\prime@\null}}
\mathchardef\prime@="0230
\let\dsize\displaystyle

\let\ssize\scriptstyle

\message{math spacing,}
\def\,{\RIfM@\mskip\thinmuskip\relax\else\kern.16667em\fi}
\def\!{\RIfM@\mskip-\thinmuskip\relax\else\kern-.16667em\fi}
\let\thinspace\,
\let\negthinspace\!
\def\medspace{\RIfM@\mskip\medmuskip\relax\else\kern.222222em\fi}
\def\negmedspace{\RIfM@\mskip-\medmuskip\relax\else\kern-.222222em\fi}
\def\thickspace{\RIfM@\mskip\thickmuskip\relax\else\kern.27777em\fi}
\let\;\thickspace
\def\negthickspace{\RIfM@\mskip-\thickmuskip\relax\else
 \kern-.27777em\fi}
\atdef@,{\RIfM@\mskip.1\thinmuskip\else\leavevmode\null,\fi}
\atdef@!{\RIfM@\mskip-.1\thinmuskip\else\leavevmode\null!\fi}
\atdef@.{\RIfM@&&\else\leavevmode.\spacefactor3000 \fi}
\def\and{\DOTSB\;\mathchar"3026 \;}

\message{fractions,}
\def\frac#1#2{{#1\over#2}}

\newdimen\ex@
\ex@.2326ex
\Invalid@\thickness
\def\thickfrac{\relaxnext@
 \DN@{\ifx\next\thickness\let\next@\nextii@\else
 \DN@{\nextii@\thickness1}\fi\next@}%
 \DNii@\thickness##1##2##3{{##2\above##1\ex@##3}}%
 \FN@\next@}

\def\thickfracwithdelims#1#2{\relaxnext@\def\ldelim@{#1}\def\rdelim@{#2}%
 \DN@{\ifx\next\thickness\let\next@\nextii@\else
 \DN@{\nextii@\thickness1}\fi\next@}%
 \DNii@\thickness##1##2##3{{##2\abovewithdelims
 \ldelim@\rdelim@##1\ex@##3}}%
 \FN@\next@}

\def\:{\nobreak\hskip.1111em\mathpunct{}\nonscript\mkern-\thinmuskip{:}\hskip
 .3333emplus.0555em\relax}
\def\snug{\unskip\kern-\mathsurround}
\message{smash commands,}
\def\topsmash{\top@true\bot@false\smash@}
\def\botsmash{\top@false\bot@true\smash@}
\newif\iftop@
\newif\ifbot@
\def\smash{\top@true\bot@true\smash@}
\def\smash@{\RIfM@\expandafter\mathpalette\expandafter\mathsm@sh\else
 \expandafter\makesm@sh\fi}
\def\finsm@sh{\iftop@\ht\z@\z@\fi\ifbot@\dp\z@\z@\fi\leavevmode\boxz@}
\message{large operator symbols,}
\def\LimitsOnSums{\global\let\slimits@\displaylimits}
\def\NoLimitsOnSums{\global\let\slimits@\nolimits}
\LimitsOnSums
\mathchardef\coprod@="1360       \def\coprod{\DOTSB\coprod@\slimits@}
\mathchardef\bigvee@="1357       \def\bigvee{\DOTSB\bigvee@\slimits@}
\mathchardef\bigwedge@="1356     \def\bigwedge{\DOTSB\bigwedge@\slimits@}
\mathchardef\biguplus@="1355     \def\biguplus{\DOTSB\biguplus@\slimits@}
\mathchardef\bigcap@="1354       \def\bigcap{\DOTSB\bigcap@\slimits@}
\mathchardef\bigcup@="1353       \def\bigcup{\DOTSB\bigcup@\slimits@}
\mathchardef\prod@="1351         \def\prod{\DOTSB\prod@\slimits@}
\mathchardef\sum@="1350          \def\sum{\DOTSB\sum@\slimits@}
\mathchardef\bigotimes@="134E    \def\bigotimes{\DOTSB\bigotimes@\slimits@}
\mathchardef\bigoplus@="134C     \def\bigoplus{\DOTSB\bigoplus@\slimits@}
\mathchardef\bigodot@="134A      \def\bigodot{\DOTSB\bigodot@\slimits@}
\mathchardef\bigsqcup@="1346     \def\bigsqcup{\DOTSB\bigsqcup@\slimits@}
\message{integrals,}
\def\LimitsOnInts{\global\let\ilimits@\displaylimits}
\def\NoLimitsOnInts{\global\let\ilimits@\nolimits}
\NoLimitsOnInts
\def\int{\DOTSI\intop\ilimits@}
\def\oint{\DOTSI\ointop\ilimits@}
\def\intic@{\mathchoice{\hskip.5em}{\hskip.4em}{\hskip.4em}{\hskip.4em}}
\def\negintic@{\mathchoice
 {\hskip-.5em}{\hskip-.4em}{\hskip-.4em}{\hskip-.4em}}
\def\intkern@{\mathchoice{\!\!\!}{\!\!}{\!\!}{\!\!}}
\def\intdots@{\mathchoice{\plaincdots@}
 {{\cdotp}\mkern1.5mu{\cdotp}\mkern1.5mu{\cdotp}}
 {{\cdotp}\mkern1mu{\cdotp}\mkern1mu{\cdotp}}
 {{\cdotp}\mkern1mu{\cdotp}\mkern1mu{\cdotp}}}
\newcount\intno@
\def\iint{\DOTSI\intno@\tw@\FN@\ints@}
\def\iiint{\DOTSI\intno@\thr@@\FN@\ints@}
\def\iiiint{\DOTSI\intno@4 \FN@\ints@}
\def\idotsint{\DOTSI\intno@\z@\FN@\ints@}
\def\ints@{\findlimits@\ints@@}
\newif\iflimtoken@
\newif\iflimits@
\def\findlimits@{\limtoken@true\ifx\next\limits\limits@true
 \else\ifx\next\nolimits\limits@false\else
 \limtoken@false\ifx\ilimits@\nolimits\limits@false\else
 \ifinner\limits@false\else\limits@true\fi\fi\fi\fi}
\def\multint@{\int\ifnum\intno@=\z@\intdots@                                
 \else\intkern@\fi                                                          
 \ifnum\intno@>\tw@\int\intkern@\fi                                         
 \ifnum\intno@>\thr@@\int\intkern@\fi                                       
 \int}                                                                      
\def\multintlimits@{\intop\ifnum\intno@=\z@\intdots@\else\intkern@\fi
 \ifnum\intno@>\tw@\intop\intkern@\fi
 \ifnum\intno@>\thr@@\intop\intkern@\fi\intop}
\def\ints@@{\iflimtoken@                                                    
 \def\ints@@@{\iflimits@\negintic@\mathop{\intic@\multintlimits@}\limits    
  \else\multint@\nolimits\fi                                                
  \eat@}                                                                    
 \else                                                                      
 \def\ints@@@{\iflimits@\negintic@
  \mathop{\intic@\multintlimits@}\limits\else
  \multint@\nolimits\fi}\fi\ints@@@}
\def\LimitsOnNames{\global\let\nlimits@\displaylimits}
\def\NoLimitsOnNames{\global\let\nlimits@\nolimits@}
\LimitsOnNames
\def\nolimits@{\relaxnext@
 \DN@{\ifx\next\limits\DN@\limits{\nolimits}\else
  \let\next@\nolimits\fi\next@}%
 \FN@\next@}
\message{operator names,}
\def\newmcodes@{\mathcode`\'"27\mathcode`\*"2A\mathcode`\."613A%
 \mathcode`\-"2D\mathcode`\/"2F\mathcode`\:"603A }
\def\operatorname#1{\mathop{\newmcodes@\kern\z@\fam\z@#1}\nolimits@}
\def\operatornamewithlimits#1{\mathop{\newmcodes@\kern\z@\fam\z@#1}\nlimits@}
\def\qopname@#1{\mathop{\fam\z@#1}\nolimits@}
\def\qopnamewl@#1{\mathop{\fam\z@#1}\nlimits@}
\def\arccos{\qopname@{arccos}}
\def\arcsin{\qopname@{arcsin}}
\def\arctan{\qopname@{arctan}}
\def\arg{\qopname@{arg}}
\def\cos{\qopname@{cos}}
\def\cosh{\qopname@{cosh}}
\def\cot{\qopname@{cot}}
\def\coth{\qopname@{coth}}
\def\csc{\qopname@{csc}}
\def\deg{\qopname@{deg}}
\def\det{\qopnamewl@{det}}
\def\dim{\qopname@{dim}}
\def\exp{\qopname@{exp}}
\def\gcd{\qopnamewl@{gcd}}
\def\hom{\qopname@{hom}}
\def\inf{\qopnamewl@{inf}}
\def\injlim{\qopnamewl@{inj\,lim}}
\def\ker{\qopname@{ker}}
\def\lg{\qopname@{lg}}
\def\lim{\qopnamewl@{lim}}
\def\liminf{\qopnamewl@{lim\,inf}}
\def\limsup{\qopnamewl@{lim\,sup}}
\def\ln{\qopname@{ln}}
\def\log{\qopname@{log}}
\def\max{\qopnamewl@{max}}
\def\min{\qopnamewl@{min}}
\def\Pr{\qopnamewl@{Pr}}
\def\projlim{\qopnamewl@{proj\,lim}}
\def\sec{\qopname@{sec}}
\def\sin{\qopname@{sin}}
\def\sinh{\qopname@{sinh}}
\def\sup{\qopnamewl@{sup}}
\def\tan{\qopname@{tan}}
\def\tanh{\qopname@{tanh}}
\def\varinjlim{\mathop{\vtop{\ialign{##\crcr
 \hfil\rm lim\hfil\crcr\noalign{\nointerlineskip}\rightarrowfill\crcr
 \noalign{\nointerlineskip\kern-\ex@}\crcr}}}}
\def\varprojlim{\mathop{\vtop{\ialign{##\crcr
 \hfil\rm lim\hfil\crcr\noalign{\nointerlineskip}\leftarrowfill\crcr
 \noalign{\nointerlineskip\kern-\ex@}\crcr}}}}
\def\varliminf{\mathop{\underline{\vrule height\z@ depth.2exwidth\z@
 \hbox{\rm lim}}}}

\newdimen\buffer@
\buffer@\fontdimen13 \tenex
\newdimen\buffer
\buffer\buffer@

\def\ResetBuffer{\fontdimen13 \tenex\buffer@\global\buffer\buffer@}
\def\shave#1{\mathop{\hbox{$\m@th\fontdimen13 \tenex\z@                     
 \displaystyle{#1}$}}\fontdimen13 \tenex\buffer}

\message{multilevel sub/superscripts,}
\Invalid@\\
\def\Let@{\relax\iffalse{\fi\let\\=\cr\iffalse}\fi}
\Invalid@\vspace
\def\vspace@{\def\vspace##1{\crcr\noalign{\vskip##1\relax}}}
\def\multilimits@{\bgroup\vspace@\Let@
 \baselineskip\fontdimen10 \scriptfont\tw@
 \advance\baselineskip\fontdimen12 \scriptfont\tw@
 \lineskip\thr@@\fontdimen8 \scriptfont\thr@@
 \lineskiplimit\lineskip
 \vbox\bgroup\ialign\bgroup\hfil$\m@th\scriptstyle{##}$\hfil\crcr}
\def\Sb{_\multilimits@}
\def\endSb{\crcr\egroup\egroup\egroup}
\def\Sp{^\multilimits@}

\def\spreadlines#1{\RIfMIfI@\onlydmatherr@\spreadlines\else
 \openup#1\relax\fi\else\onlydmatherr@\spreadlines\fi}
\def\Mathstrut@{\copy\Mathstrutbox@}
\newbox\Mathstrutbox@
\setbox\Mathstrutbox@\null
\setboxz@h{$\m@th($}
\ht\Mathstrutbox@\ht\z@
\dp\Mathstrutbox@\dp\z@
\message{matrices,}
\newdimen\spreadmlines@
\def\spreadmatrixlines#1{\RIfMIfI@
 \onlydmatherr@\spreadmatrixlines\else
 \spreadmlines@#1\relax\fi\else\onlydmatherr@\spreadmatrixlines\fi}
\def\matrix{\null\,\vcenter\bgroup\Let@\vspace@
 \normalbaselines\openup\spreadmlines@\ialign
 \bgroup\hfil$\m@th##$\hfil&&\quad\hfil$\m@th##$\hfil\crcr
 \Mathstrut@\crcr\noalign{\kern-\baselineskip}}
\def\endmatrix{\crcr\Mathstrut@\crcr\noalign{\kern-\baselineskip}\egroup
 \egroup\,}
\def\format{\crcr\egroup\iffalse{\fi\ifnum`}=0 \fi\format@}
\newtoks\hashtoks@
\hashtoks@{#}
\def\format@#1\\{\def\preamble@{#1}%
 \def\l{$\m@th\the\hashtoks@$\hfil}%
 \def\c{\hfil$\m@th\the\hashtoks@$\hfil}%
 \def\r{\hfil$\m@th\the\hashtoks@$}%
 \edef\preamble@@{\preamble@}\ifnum`{=0 \fi\iffalse}\fi
 \ialign\bgroup\span\preamble@@\crcr}
\def\smallmatrix{\null\,\vcenter\bgroup\vspace@\Let@
 \baselineskip9\ex@\lineskip\ex@
 \ialign\bgroup\hfil$\m@th\scriptstyle{##}$\hfil&&\thickspace\hfil
 $\m@th\scriptstyle{##}$\hfil\crcr}
\def\endsmallmatrix{\crcr\egroup\egroup\,}

\newmuskip\dotsspace@
\dotsspace@1.5mu
\def\strip@#1 {#1}
\def\spacehdots#1\for#2{\multispan{#2}\xleaders
 \hbox{$\m@th\mkern\strip@#1 \dotsspace@.\mkern\strip@#1 \dotsspace@$}\hfill}
\def\hdotsfor#1{\spacehdots\@ne\for{#1}}
\def\multispan@#1{\omit\mscount#1\unskip\loop\ifnum\mscount>\@ne\sp@n\repeat}
\def\spaceinnerhdots#1\for#2\after#3{\multispan@{\strip@#2 }#3\xleaders
 \hbox{$\m@th\mkern\strip@#1 \dotsspace@.\mkern\strip@#1 \dotsspace@$}\hfill}
\def\innerhdotsfor#1\after#2{\spaceinnerhdots\@ne\for#1\after{#2}}
\def\cases{\bgroup\spreadmlines@\jot\left\{\,\matrix\format\l&\quad\l\\}
\def\endcases{\endmatrix\right.\egroup}
\message{multiline displays,}
\newif\ifinany@
\newif\ifinalign@
\newif\ifingather@
\def\strut@{\copy\strutbox@}
\newbox\strutbox@
\setbox\strutbox@\hbox{\vrule height8\p@ depth3\p@ width\z@}
\def\topaligned{\null\,\vtop\aligned@}
\def\botaligned{\null\,\vbox\aligned@}
\def\aligned{\null\,\vcenter\aligned@}
\def\aligned@{\bgroup\vspace@\Let@
 \ifinany@\else\openup\jot\fi\ialign
 \bgroup\hfil\strut@$\m@th\displaystyle{##}$&
 $\m@th\displaystyle{{}##}$\hfil\crcr}
\def\endaligned{\crcr\egroup\egroup}

\def\alignedat#1{\null\,\vcenter\bgroup\doat@{#1}\vspace@\Let@
 \ifinany@\else\openup\jot\fi\ialign\bgroup\span\preamble@@\crcr}
\newcount\atcount@
\def\doat@#1{\toks@{\hfil\strut@$\m@th
 \displaystyle{\the\hashtoks@}$&$\m@th\displaystyle
 {{}\the\hashtoks@}$\hfil}
 \atcount@#1\relax\advance\atcount@\m@ne                                    
 \loop\ifnum\atcount@>\z@\toks@=\expandafter{\the\toks@&\hfil$\m@th
 \displaystyle{\the\hashtoks@}$&$\m@th
 \displaystyle{{}\the\hashtoks@}$\hfil}\advance
  \atcount@\m@ne\repeat                                                     
 \xdef\preamble@{\the\toks@}\xdef\preamble@@{\preamble@}}

\def\gathered{\null\,\vcenter\bgroup\vspace@\Let@
 \ifinany@\else\openup\jot\fi\ialign
 \bgroup\hfil\strut@$\m@th\displaystyle{##}$\hfil\crcr}
\def\endgathered{\crcr\egroup\egroup}
\newif\iftagsleft@
\def\TagsOnLeft{\global\tagsleft@true}
\def\TagsOnRight{\global\tagsleft@false}
\TagsOnLeft
\newif\ifmathtags@
\def\TagsAsMath{\global\mathtags@true}
\def\TagsAsText{\global\mathtags@false}
\TagsAsText
\def\tagform@#1{\hbox{\rm(\ignorespaces#1\unskip)}}
\def\thetag{\leavevmode\tagform@}
\def\tag#1$${\iftagsleft@\leqno\else\eqno\fi                                
 \maketag@#1\maketag@                                                       
 $$}                                                                        
\def\maketag@{\FN@\maketag@@}
\def\maketag@@{\ifx\next"\expandafter\maketag@@@\else\expandafter\maketag@@@@
 \fi}
\def\maketag@@@"#1"#2\maketag@{\hbox{\rm#1}}                                
\def\maketag@@@@#1\maketag@{\ifmathtags@\tagform@{$\m@th#1$}\else
 \tagform@{#1}\fi}
\interdisplaylinepenalty\@M
\def\allowdisplaybreaks{\RIfMIfI@
 \onlydmatherr@\allowdisplaybreaks\else
 \interdisplaylinepenalty\z@\fi\else\onlydmatherr@\allowdisplaybreaks\fi}
\Invalid@\allowdisplaybreak
\Invalid@\displaybreak
\Invalid@\intertext
\def\allowdisplaybreak@{\def\allowdisplaybreak{\crcr\noalign{\allowbreak}}}
\def\displaybreak@{\def\displaybreak{\crcr\noalign{\break}}}
\def\intertext@{\def\intertext##1{\crcr\noalign{%
 \penalty\postdisplaypenalty \vskip\belowdisplayskip
 \vbox{\normalbaselines\noindent##1}%
 \penalty\predisplaypenalty \vskip\abovedisplayskip}}}
\newskip\centering@
\centering@\z@ plus\@m\p@
\def\align{\relax\ifingather@\DN@{\csname align (in
  \string\gather)\endcsname}\else
 \ifmmode\ifinner\DN@{\onlydmatherr@\align}\else
  \let\next@\align@\fi
 \else\DN@{\onlydmatherr@\align}\fi\fi\next@}
\newhelp\andhelp@
{An extra & here is so disastrous that you should probably exit^^J
and fix things up.}
\newif\iftag@
\newcount\and@
\def\align@{\inalign@true\inany@true
 \vspace@\allowdisplaybreak@\displaybreak@\intertext@
 \def\tag{\global\tag@true\ifnum\and@=\z@\DN@{&&}\else
          \DN@{&}\fi\next@}%
 \iftagsleft@\DN@{\csname align \endcsname}\else
  \DN@{\csname align \space\endcsname}\fi\next@}
\def\Tag@{\iftag@\else\errhelp\andhelp@\err@{Extra & on this line}\fi}
\newdimen\lwidth@
\newdimen\rwidth@
\newdimen\maxlwidth@
\newdimen\maxrwidth@
\newdimen\totwidth@
\def\measure@#1\endalign{\lwidth@\z@\rwidth@\z@\maxlwidth@\z@\maxrwidth@\z@
 \global\and@\z@                                                            
 \setbox@ne\vbox                                                            
  {\everycr{\noalign{\global\tag@false\global\and@\z@}}\Let@                
  \halign{\setboxz@h{$\m@th\displaystyle{\@lign##}$}
   \global\lwidth@\wdz@                                                     
   \ifdim\lwidth@>\maxlwidth@\global\maxlwidth@\lwidth@\fi                  
   \global\advance\and@\@ne                                                 
   &\setboxz@h{$\m@th\displaystyle{{}\@lign##}$}\global\rwidth@\wdz@        
   \ifdim\rwidth@>\maxrwidth@\global\maxrwidth@\rwidth@\fi                  
   \global\advance\and@\@ne                                                
   &\Tag@
   \eat@{##}\crcr#1\crcr}}
 \totwidth@\maxlwidth@\advance\totwidth@\maxrwidth@}                       
\def\displ@y@{\global\dt@ptrue\openup\jot
 \everycr{\noalign{\global\tag@false\global\and@\z@\ifdt@p\global\dt@pfalse
 \vskip-\lineskiplimit\vskip\normallineskiplimit\else
 \penalty\interdisplaylinepenalty\fi}}}
\def\black@#1{\noalign{\ifdim#1>\displaywidth
 \dimen@\prevdepth\nointerlineskip                                          
 \vskip-\ht\strutbox@\vskip-\dp\strutbox@                                   
 \vbox{\noindent\hbox to#1{\strut@\hfill}}
 \prevdepth\dimen@                                                          
 \fi}}
\expandafter\def\csname align \space\endcsname#1\endalign
 {\measure@#1\endalign\global\and@\z@                                       
 \ifingather@\everycr{\noalign{\global\and@\z@}}\else\displ@y@\fi           
 \Let@\tabskip\centering@                                                   
 \halign to\displaywidth
  {\hfil\strut@\setboxz@h{$\m@th\displaystyle{\@lign##}$}
  \global\lwidth@\wdz@\boxz@\global\advance\and@\@ne                        
  \tabskip\z@skip                                                           
  &\setboxz@h{$\m@th\displaystyle{{}\@lign##}$}
  \global\rwidth@\wdz@\boxz@\hfill\global\advance\and@\@ne                  
  \tabskip\centering@                                                       
  &\setboxz@h{\@lign\strut@\maketag@##\maketag@}
  \dimen@\displaywidth\advance\dimen@-\totwidth@
  \divide\dimen@\tw@\advance\dimen@\maxrwidth@\advance\dimen@-\rwidth@     
  \ifdim\dimen@<\tw@\wdz@\llap{\vtop{\normalbaselines\null\boxz@}}
  \else\llap{\boxz@}\fi                                                    
  \tabskip\z@skip                                                          
  \crcr#1\crcr                                                             
  \black@\totwidth@}}                                                      
\newdimen\lineht@
\expandafter\def\csname align \endcsname#1\endalign{\measure@#1\endalign
 \global\and@\z@
 \ifdim\totwidth@>\displaywidth\let\displaywidth@\totwidth@\else
  \let\displaywidth@\displaywidth\fi                                        
 \ifingather@\everycr{\noalign{\global\and@\z@}}\else\displ@y@\fi
 \Let@\tabskip\centering@\halign to\displaywidth
  {\hfil\strut@\setboxz@h{$\m@th\displaystyle{\@lign##}$}%
  \global\lwidth@\wdz@\global\lineht@\ht\z@                                 
  \boxz@\global\advance\and@\@ne
  \tabskip\z@skip&\setboxz@h{$\m@th\displaystyle{{}\@lign##}$}%
  \global\rwidth@\wdz@\ifdim\ht\z@>\lineht@\global\lineht@\ht\z@\fi         
  \boxz@\hfil\global\advance\and@\@ne
  \tabskip\centering@&\kern-\displaywidth@                                  
  \setboxz@h{\@lign\strut@\maketag@##\maketag@}%
  \dimen@\displaywidth\advance\dimen@-\totwidth@
  \divide\dimen@\tw@\advance\dimen@\maxlwidth@\advance\dimen@-\lwidth@
  \ifdim\dimen@<\tw@\wdz@
   \rlap{\vbox{\normalbaselines\boxz@\vbox to\lineht@{}}}\else
   \rlap{\boxz@}\fi
  \tabskip\displaywidth@\crcr#1\crcr\black@\totwidth@}}
\expandafter\def\csname align (in \string\gather)\endcsname
  #1\endalign{\vcenter{\align@#1\endalign}}
\Invalid@\endalign
\newif\ifxat@
\def\alignat{\RIfMIfI@\DN@{\onlydmatherr@\alignat}\else
 \DN@{\csname alignat \endcsname}\fi\else
 \DN@{\onlydmatherr@\alignat}\fi\next@}
\newif\ifmeasuring@
\newbox\savealignat@
\expandafter\def\csname alignat \endcsname#1#2\endalignat                   
 {\inany@true\xat@false
 \def\tag{\global\tag@true\count@#1\relax\multiply\count@\tw@
  \xdef\tag@{}\loop\ifnum\count@>\and@\xdef\tag@{&\tag@}\advance\count@\m@ne
  \repeat\tag@}%
 \vspace@\allowdisplaybreak@\displaybreak@\intertext@
 \displ@y@\measuring@true                                                   
 \setbox\savealignat@\hbox{$\m@th\displaystyle\Let@
  \attag@{#1}
  \vbox{\halign{\span\preamble@@\crcr#2\crcr}}$}%
 \measuring@false                                                           
 \Let@\attag@{#1}
 \tabskip\centering@\halign to\displaywidth
  {\span\preamble@@\crcr#2\crcr                                             
  \black@{\wd\savealignat@}}}                                               
\Invalid@\endalignat
\def\xalignat{\RIfMIfI@
 \DN@{\onlydmatherr@\xalignat}\else
 \DN@{\csname xalignat \endcsname}\fi\else
 \DN@{\onlydmatherr@\xalignat}\fi\next@}
\expandafter\def\csname xalignat \endcsname#1#2\endxalignat
 {\inany@true\xat@true
 \def\tag{\global\tag@true\def\tag@{}\count@#1\relax\multiply\count@\tw@
  \loop\ifnum\count@>\and@\xdef\tag@{&\tag@}\advance\count@\m@ne\repeat\tag@}%
 \vspace@\allowdisplaybreak@\displaybreak@\intertext@
 \displ@y@\measuring@true\setbox\savealignat@\hbox{$\m@th\displaystyle\Let@
 \attag@{#1}\vbox{\halign{\span\preamble@@\crcr#2\crcr}}$}%
 \measuring@false\Let@
 \attag@{#1}\tabskip\centering@\halign to\displaywidth
 {\span\preamble@@\crcr#2\crcr\black@{\wd\savealignat@}}}
\def\attag@#1{\let\Maketag@\maketag@\let\TAG@\Tag@                          
 \let\Tag@=0\let\maketag@=0
 \ifmeasuring@\def\llap@##1{\setboxz@h{##1}\hbox to\tw@\wdz@{}}%
  \def\rlap@##1{\setboxz@h{##1}\hbox to\tw@\wdz@{}}\else
  \let\llap@\llap\let\rlap@\rlap\fi                                         
 \toks@{\hfil\strut@$\m@th\displaystyle{\@lign\the\hashtoks@}$\tabskip\z@skip
  \global\advance\and@\@ne&$\m@th\displaystyle{{}\@lign\the\hashtoks@}$\hfil
  \ifxat@\tabskip\centering@\fi\global\advance\and@\@ne}
 \iftagsleft@
  \toks@@{\tabskip\centering@&\Tag@\kern-\displaywidth
   \rlap@{\@lign\maketag@\the\hashtoks@\maketag@}%
   \global\advance\and@\@ne\tabskip\displaywidth}\else
  \toks@@{\tabskip\centering@&\Tag@\llap@{\@lign\maketag@
   \the\hashtoks@\maketag@}\global\advance\and@\@ne\tabskip\z@skip}\fi      
 \atcount@#1\relax\advance\atcount@\m@ne
 \loop\ifnum\atcount@>\z@
 \toks@=\expandafter{\the\toks@&\hfil$\m@th\displaystyle{\@lign
  \the\hashtoks@}$\global\advance\and@\@ne
  \tabskip\z@skip&$\m@th\displaystyle{{}\@lign\the\hashtoks@}$\hfil\ifxat@
  \tabskip\centering@\fi\global\advance\and@\@ne}\advance\atcount@\m@ne
 \repeat                                                                    
 \xdef\preamble@{\the\toks@\the\toks@@}
 \xdef\preamble@@{\preamble@}
 \let\maketag@\Maketag@\let\Tag@\TAG@}                                      
\Invalid@\endxalignat
\def\xxalignat{\RIfMIfI@
 \DN@{\onlydmatherr@\xxalignat}\else\DN@{\csname xxalignat
  \endcsname}\fi\else
 \DN@{\onlydmatherr@\xxalignat}\fi\next@}
\expandafter\def\csname xxalignat \endcsname#1#2\endxxalignat{\inany@true
 \vspace@\allowdisplaybreak@\displaybreak@\intertext@
 \displ@y\setbox\savealignat@\hbox{$\m@th\displaystyle\Let@
 \xxattag@{#1}\vbox{\halign{\span\preamble@@\crcr#2\crcr}}$}%
 \Let@\xxattag@{#1}\tabskip\z@skip\halign to\displaywidth
 {\span\preamble@@\crcr#2\crcr\black@{\wd\savealignat@}}}
\def\xxattag@#1{\toks@{\tabskip\z@skip\hfil\strut@
 $\m@th\displaystyle{\the\hashtoks@}$&%
 $\m@th\displaystyle{{}\the\hashtoks@}$\hfil\tabskip\centering@&}%
 \atcount@#1\relax\advance\atcount@\m@ne\loop\ifnum\atcount@>\z@
 \toks@=\expandafter{\the\toks@&\hfil$\m@th\displaystyle{\the\hashtoks@}$%
  \tabskip\z@skip&$\m@th\displaystyle{{}\the\hashtoks@}$\hfil
  \tabskip\centering@}\advance\atcount@\m@ne\repeat
 \xdef\preamble@{\the\toks@\tabskip\z@skip}\xdef\preamble@@{\preamble@}}
\Invalid@\endxxalignat
\newdimen\gwidth@
\newdimen\gmaxwidth@
\def\gmeasure@#1\endgather{\gwidth@\z@\gmaxwidth@\z@\setbox@ne\vbox{\Let@
 \halign{\setboxz@h{$\m@th\displaystyle{##}$}\global\gwidth@\wdz@
 \ifdim\gwidth@>\gmaxwidth@\global\gmaxwidth@\gwidth@\fi
 &\eat@{##}\crcr#1\crcr}}}
\def\gather{\RIfMIfI@\DN@{\onlydmatherr@\gather}\else
 \ingather@true\inany@true\def\tag{&}%
 \vspace@\allowdisplaybreak@\displaybreak@\intertext@
 \displ@y\Let@
 \iftagsleft@\DN@{\csname gather \endcsname}\else
  \DN@{\csname gather \space\endcsname}\fi\fi
 \else\DN@{\onlydmatherr@\gather}\fi\next@}
\expandafter\def\csname gather \space\endcsname#1\endgather
 {\gmeasure@#1\endgather\tabskip\centering@
 \halign to\displaywidth{\hfil\strut@\setboxz@h{$\m@th\displaystyle{##}$}%
 \global\gwidth@\wdz@\boxz@\hfil&
 \setboxz@h{\strut@{\maketag@##\maketag@}}%
 \dimen@\displaywidth\advance\dimen@-\gwidth@
 \ifdim\dimen@>\tw@\wdz@\llap{\boxz@}\else
 \llap{\vtop{\normalbaselines\null\boxz@}}\fi
 \tabskip\z@skip\crcr#1\crcr\black@\gmaxwidth@}}
\newdimen\glineht@
\expandafter\def\csname gather \endcsname#1\endgather{\gmeasure@#1\endgather
 \ifdim\gmaxwidth@>\displaywidth\let\gdisplaywidth@\gmaxwidth@\else
 \let\gdisplaywidth@\displaywidth\fi\tabskip\centering@\halign to\displaywidth
 {\hfil\strut@\setboxz@h{$\m@th\displaystyle{##}$}%
 \global\gwidth@\wdz@\global\glineht@\ht\z@\boxz@\hfil&\kern-\gdisplaywidth@
 \setboxz@h{\strut@{\maketag@##\maketag@}}%
 \dimen@\displaywidth\advance\dimen@-\gwidth@
 \ifdim\dimen@>\tw@\wdz@\rlap{\boxz@}\else
 \rlap{\vbox{\normalbaselines\boxz@\vbox to\glineht@{}}}\fi
 \tabskip\gdisplaywidth@\crcr#1\crcr\black@\gmaxwidth@}}
\newif\ifctagsplit@
\def\CenteredTagsOnSplits{\global\ctagsplit@true}
\def\TopOrBottomTagsOnSplits{\global\ctagsplit@false}
\TopOrBottomTagsOnSplits
\def\split{\relax\ifinany@\let\next@\insplit@\else
 \ifmmode\ifinner\def\next@{\onlydmatherr@\split}\else
 \let\next@\outsplit@\fi\else
 \def\next@{\onlydmatherr@\split}\fi\fi\next@}
\def\insplit@{\global\setbox\z@\vbox\bgroup\vspace@\Let@\ialign\bgroup
 \hfil\strut@$\m@th\displaystyle{##}$&$\m@th\displaystyle{{}##}$\hfill\crcr}
\def\endsplit{\crcr\egroup\egroup\iftagsleft@\expandafter\lendsplit@\else
 \expandafter\rendsplit@\fi}
\def\rendsplit@{\global\setbox9 \vbox
 {\unvcopy\z@\global\setbox8 \lastbox\unskip}
 \setbox@ne\hbox{\unhcopy8 \unskip\global\setbox\tw@\lastbox
 \unskip\global\setbox\thr@@\lastbox}
 \global\setbox7 \hbox{\unhbox\tw@\unskip}
 \ifinalign@\ifctagsplit@                                                   
  \gdef\split@{\hbox to\wd\thr@@{}&
   \vcenter{\vbox{\moveleft\wd\thr@@\boxz@}}}
 \else\gdef\split@{&\vbox{\moveleft\wd\thr@@\box9}\crcr
  \box\thr@@&\box7}\fi                                                      
 \else                                                                      
  \ifctagsplit@\gdef\split@{\vcenter{\boxz@}}\else
  \gdef\split@{\box9\crcr\hbox{\box\thr@@\box7}}\fi
 \fi
 \split@}                                                                   
\def\lendsplit@{\global\setbox9\vtop{\unvcopy\z@}
 \setbox@ne\vbox{\unvcopy\z@\global\setbox8\lastbox}
 \setbox@ne\hbox{\unhcopy8\unskip\setbox\tw@\lastbox
  \unskip\global\setbox\thr@@\lastbox}
 \ifinalign@\ifctagsplit@                                                   
  \gdef\split@{\hbox to\wd\thr@@{}&
  \vcenter{\vbox{\moveleft\wd\thr@@\box9}}}
  \else                                                                     
  \gdef\split@{\hbox to\wd\thr@@{}&\vbox{\moveleft\wd\thr@@\box9}}\fi
 \else
  \ifctagsplit@\gdef\split@{\vcenter{\box9}}\else
  \gdef\split@{\box9}\fi
 \fi\split@}
\def\outsplit@#1$${\align\insplit@#1\endalign$$}
\newdimen\multlinegap@
\multlinegap@1em
\newdimen\multlinetaggap@
\multlinetaggap@1em
\def\MultlineGap#1{\global\multlinegap@#1\relax}
\def\multlinegap#1{\RIfMIfI@\onlydmatherr@\multlinegap\else
 \multlinegap@#1\relax\fi\else\onlydmatherr@\multlinegap\fi}
\def\nomultlinegap{\multlinegap{\z@}}
\def\multline{\RIfMIfI@
 \DN@{\onlydmatherr@\multline}\else
 \DN@{\multline@}\fi\else
 \DN@{\onlydmatherr@\multline}\fi\next@}
\newif\iftagin@
\def\tagin@#1{\tagin@false\in@\tag{#1}\ifin@\tagin@true\fi}
\def\multline@#1$${\inany@true\vspace@\allowdisplaybreak@\displaybreak@
 \tagin@{#1}\iftagsleft@\DN@{\multline@l#1$$}\else
 \DN@{\multline@r#1$$}\fi\next@}
\newdimen\mwidth@
\def\rmmeasure@#1\endmultline{%
 \def\shoveleft##1{##1}\def\shoveright##1{##1}
 \setbox@ne\vbox{\Let@\halign{\setboxz@h
  {$\m@th\@lign\displaystyle{}##$}\global\mwidth@\wdz@
  \crcr#1\crcr}}}
\newdimen\mlineht@
\newif\ifzerocr@
\newif\ifonecr@
\def\lmmeasure@#1\endmultline{\global\zerocr@true\global\onecr@false
 \everycr{\noalign{\ifonecr@\global\onecr@false\fi
  \ifzerocr@\global\zerocr@false\global\onecr@true\fi}}
  \def\shoveleft##1{##1}\def\shoveright##1{##1}%
 \setbox@ne\vbox{\Let@\halign{\setboxz@h
  {$\m@th\@lign\displaystyle{}##$}\ifonecr@\global\mwidth@\wdz@
  \global\mlineht@\ht\z@\fi\crcr#1\crcr}}}
\newbox\mtagbox@
\newdimen\ltwidth@
\newdimen\rtwidth@
\def\multline@l#1$${\iftagin@\DN@{\lmultline@@#1$$}\else
 \DN@{\setbox\mtagbox@\null\ltwidth@\z@\rtwidth@\z@
  \lmultline@@@#1$$}\fi\next@}
\def\lmultline@@#1\endmultline\tag#2$${%
 \setbox\mtagbox@\hbox{\maketag@#2\maketag@}
 \lmmeasure@#1\endmultline\dimen@\mwidth@\advance\dimen@\wd\mtagbox@
 \advance\dimen@\multlinetaggap@                                            
 \ifdim\dimen@>\displaywidth\ltwidth@\z@\else\ltwidth@\wd\mtagbox@\fi       
 \lmultline@@@#1\endmultline$$}
\def\lmultline@@@{\displ@y
 \def\shoveright##1{##1\hfilneg\hskip\multlinegap@}%
 \def\shoveleft##1{\setboxz@h{$\m@th\displaystyle{}##1$}%
  \setbox@ne\hbox{$\m@th\displaystyle##1$}%
  \hfilneg
  \iftagin@
   \ifdim\ltwidth@>\z@\hskip\ltwidth@\hskip\multlinetaggap@\fi
  \else\hskip\multlinegap@\fi\hskip.5\wd@ne\hskip-.5\wdz@##1}
  \halign\bgroup\Let@\hbox to\displaywidth
   {\strut@$\m@th\displaystyle\hfil{}##\hfil$}\crcr
   \hfilneg                                                                 
   \iftagin@                                                                
    \ifdim\ltwidth@>\z@                                                     
     \box\mtagbox@\hskip\multlinetaggap@                                    
    \else
     \rlap{\vbox{\normalbaselines\hbox{\strut@\box\mtagbox@}%
     \vbox to\mlineht@{}}}\fi                                               
   \else\hskip\multlinegap@\fi}                                             
\def\multline@r#1$${\iftagin@\DN@{\rmultline@@#1$$}\else
 \DN@{\setbox\mtagbox@\null\ltwidth@\z@\rtwidth@\z@
  \rmultline@@@#1$$}\fi\next@}
\def\rmultline@@#1\endmultline\tag#2$${\ltwidth@\z@
 \setbox\mtagbox@\hbox{\maketag@#2\maketag@}%
 \rmmeasure@#1\endmultline\dimen@\mwidth@\advance\dimen@\wd\mtagbox@
 \advance\dimen@\multlinetaggap@
 \ifdim\dimen@>\displaywidth\rtwidth@\z@\else\rtwidth@\wd\mtagbox@\fi
 \rmultline@@@#1\endmultline$$}
\def\rmultline@@@{\displ@y
 \def\shoveright##1{##1\hfilneg\iftagin@\ifdim\rtwidth@>\z@
  \hskip\rtwidth@\hskip\multlinetaggap@\fi\else\hskip\multlinegap@\fi}%
 \def\shoveleft##1{\setboxz@h{$\m@th\displaystyle{}##1$}%
  \setbox@ne\hbox{$\m@th\displaystyle##1$}%
  \hfilneg\hskip\multlinegap@\hskip.5\wd@ne\hskip-.5\wdz@##1}%
 \halign\bgroup\Let@\hbox to\displaywidth
  {\strut@$\m@th\displaystyle\hfil{}##\hfil$}\crcr
 \hfilneg\hskip\multlinegap@}
\def\endmultline{\iftagsleft@\expandafter\lendmultline@\else
 \expandafter\rendmultline@\fi}
\def\lendmultline@{\hfilneg\hskip\multlinegap@\crcr\egroup}
\def\rendmultline@{\iftagin@                                                
 \ifdim\rtwidth@>\z@                                                        
  \hskip\multlinetaggap@\box\mtagbox@                                       
 \else\llap{\vtop{\normalbaselines\null\hbox{\strut@\box\mtagbox@}}}\fi     
 \else\hskip\multlinegap@\fi                                                
 \hfilneg\crcr\egroup}
\def\bmod{\mskip-\medmuskip\mkern5mu\mathbin{\fam\z@ mod}\penalty900
 \mkern5mu\mskip-\medmuskip}
\def\pmod#1{\allowbreak\ifinner\mkern8mu\else\mkern18mu\fi
 ({\fam\z@ mod}\,\,#1)}
\def\pod#1{\allowbreak\ifinner\mkern8mu\else\mkern18mu\fi(#1)}
\def\mod#1{\allowbreak\ifinner\mkern12mu\else\mkern18mu\fi{\fam\z@ mod}\,\,#1}
\message{continued fractions,}
\newcount\cfraccount@
\def\cfrac{\bgroup\bgroup\advance\cfraccount@\@ne\strut
 \iffalse{\fi\def\\{\over\displaystyle}\iffalse}\fi}
\def\lcfrac{\bgroup\bgroup\advance\cfraccount@\@ne\strut
 \iffalse{\fi\def\\{\hfill\over\displaystyle}\iffalse}\fi}
\def\rcfrac{\bgroup\bgroup\advance\cfraccount@\@ne\strut\hfill
 \iffalse{\fi\def\\{\over\displaystyle}\iffalse}\fi}
\def\gloop@#1\repeat{\gdef\body{#1}\iterate}
\def\endcfrac{\gloop@\ifnum\cfraccount@>\z@\global\advance\cfraccount@\m@ne
 \egroup\hskip-\nulldelimiterspace\egroup\repeat}
\message{compound symbols,}
\def\binrel@#1{\setboxz@h{\thinmuskip0mu
  \medmuskip\m@ne mu\thickmuskip\@ne mu$#1\m@th$}%
 \setbox@ne\hbox{\thinmuskip0mu\medmuskip\m@ne mu\thickmuskip
  \@ne mu${}#1{}\m@th$}%
 \setbox\tw@\hbox{\hskip\wd@ne\hskip-\wdz@}}
\def\overset#1\to#2{\binrel@{#2}\ifdim\wd\tw@<\z@
 \mathbin{\mathop{\kern\z@#2}\limits^{#1}}\else\ifdim\wd\tw@>\z@
 \mathrel{\mathop{\kern\z@#2}\limits^{#1}}\else
 {\mathop{\kern\z@#2}\limits^{#1}}{}\fi\fi}
\def\underset#1\to#2{\binrel@{#2}\ifdim\wd\tw@<\z@
 \mathbin{\mathop{\kern\z@#2}\limits_{#1}}\else\ifdim\wd\tw@>\z@
 \mathrel{\mathop{\kern\z@#2}\limits_{#1}}\else
 {\mathop{\kern\z@#2}\limits_{#1}}{}\fi\fi}
\def\oversetbrace#1\to#2{\overbrace{#2}^{#1}}
\def\undersetbrace#1\to#2{\underbrace{#2}_{#1}}
\def\sideset#1\and#2\to#3{%
 \setbox@ne\hbox{$\dsize{\vphantom{#3}}#1{#3}\m@th$}%
 \setbox\tw@\hbox{$\dsize{#3}#2\m@th$}%
 \hskip\wd@ne\hskip-\wd\tw@\mathop{\hskip\wd\tw@\hskip-\wd@ne
  {\vphantom{#3}}#1{#3}#2}}
\def\rightarrowfill@#1{\setboxz@h{$#1-\m@th$}\ht\z@\z@
  $#1\m@th\copy\z@\mkern-6mu\cleaders
  \hbox{$#1\mkern-2mu\box\z@\mkern-2mu$}\hfill
  \mkern-6mu\mathord\rightarrow$}
\def\leftarrowfill@#1{\setboxz@h{$#1-\m@th$}\ht\z@\z@
  $#1\m@th\mathord\leftarrow\mkern-6mu\cleaders
  \hbox{$#1\mkern-2mu\copy\z@\mkern-2mu$}\hfill
  \mkern-6mu\box\z@$}
\def\leftrightarrowfill@#1{\setboxz@h{$#1-\m@th$}\ht\z@\z@
  $#1\m@th\mathord\leftarrow\mkern-6mu\cleaders
  \hbox{$#1\mkern-2mu\box\z@\mkern-2mu$}\hfill
  \mkern-6mu\mathord\rightarrow$}
\def\overrightarrow{\mathpalette\overrightarrow@}
\def\overrightarrow@#1#2{\vbox{\ialign{##\crcr\rightarrowfill@#1\crcr
 \noalign{\kern-\ex@\nointerlineskip}$\m@th\hfil#1#2\hfil$\crcr}}}

\def\overleftarrow{\mathpalette\overleftarrow@}
\def\overleftarrow@#1#2{\vbox{\ialign{##\crcr\leftarrowfill@#1\crcr
 \noalign{\kern-\ex@\nointerlineskip}$\m@th\hfil#1#2\hfil$\crcr}}}
\def\overleftrightarrow{\mathpalette\overleftrightarrow@}
\def\overleftrightarrow@#1#2{\vbox{\ialign{##\crcr\leftrightarrowfill@#1\crcr
 \noalign{\kern-\ex@\nointerlineskip}$\m@th\hfil#1#2\hfil$\crcr}}}
\def\underrightarrow{\mathpalette\underrightarrow@}
\def\underrightarrow@#1#2{\vtop{\ialign{##\crcr$\m@th\hfil#1#2\hfil$\crcr
 \noalign{\nointerlineskip}\rightarrowfill@#1\crcr}}}

\def\underleftarrow{\mathpalette\underleftarrow@}
\def\underleftarrow@#1#2{\vtop{\ialign{##\crcr$\m@th\hfil#1#2\hfil$\crcr
 \noalign{\nointerlineskip}\leftarrowfill@#1\crcr}}}
\def\underleftrightarrow{\mathpalette\underleftrightarrow@}
\def\underleftrightarrow@#1#2{\vtop{\ialign{##\crcr$\m@th\hfil#1#2\hfil$\crcr
 \noalign{\nointerlineskip}\leftrightarrowfill@#1\crcr}}}
\message{various kinds of dots,}
\let\DOTSI\relax
\let\DOTSB\relax

\newif\ifmath@
{\uccode`7=`\\ \uccode`8=`m \uccode`9=`a \uccode`0=`t \uccode`!=`h
 \uppercase{\gdef\math@#1#2#3#4#5#6\math@{\global\math@false\ifx 7#1\ifx 8#2%
 \ifx 9#3\ifx 0#4\ifx !#5\xdef\meaning@{#6}\global\math@true\fi\fi\fi\fi\fi}}}
\newif\ifmathch@
{\uccode`7=`c \uccode`8=`h \uccode`9=`\"
 \uppercase{\gdef\mathch@#1#2#3#4#5#6\mathch@{\global\mathch@false
  \ifx 7#1\ifx 8#2\ifx 9#5\global\mathch@true\xdef\meaning@{9#6}\fi\fi\fi}}}
\newcount\classnum@
\def\getmathch@#1.#2\getmathch@{\classnum@#1 \divide\classnum@4096
 \ifcase\number\classnum@\or\or\gdef\thedots@{\dotsb@}\or
 \gdef\thedots@{\dotsb@}\fi}
\newif\ifmathbin@
{\uccode`4=`b \uccode`5=`i \uccode`6=`n
 \uppercase{\gdef\mathbin@#1#2#3{\relaxnext@
  \DNii@##1\mathbin@{\ifx\space@\next\global\mathbin@true\fi}%
 \global\mathbin@false\DN@##1\mathbin@{}%
 \ifx 4#1\ifx 5#2\ifx 6#3\DN@{\FN@\nextii@}\fi\fi\fi\next@}}}
\newif\ifmathrel@
{\uccode`4=`r \uccode`5=`e \uccode`6=`l
 \uppercase{\gdef\mathrel@#1#2#3{\relaxnext@
  \DNii@##1\mathrel@{\ifx\space@\next\global\mathrel@true\fi}%
 \global\mathrel@false\DN@##1\mathrel@{}%
 \ifx 4#1\ifx 5#2\ifx 6#3\DN@{\FN@\nextii@}\fi\fi\fi\next@}}}
\newif\ifmacro@
{\uccode`5=`m \uccode`6=`a \uccode`7=`c
 \uppercase{\gdef\macro@#1#2#3#4\macro@{\global\macro@false
  \ifx 5#1\ifx 6#2\ifx 7#3\global\macro@true
  \xdef\meaning@{\macro@@#4\macro@@}\fi\fi\fi}}}
\def\macro@@#1->#2\macro@@{#2}
\newif\ifDOTS@
\newcount\DOTSCASE@
{\uccode`6=`\\ \uccode`7=`D \uccode`8=`O \uccode`9=`T \uccode`0=`S
 \uppercase{\gdef\DOTS@#1#2#3#4#5{\global\DOTS@false\DN@##1\DOTS@{}%
  \ifx 6#1\ifx 7#2\ifx 8#3\ifx 9#4\ifx 0#5\let\next@\DOTS@@\fi\fi\fi\fi\fi
  \next@}}}
{\uccode`3=`B \uccode`4=`I \uccode`5=`X
 \uppercase{\gdef\DOTS@@#1{\relaxnext@
  \DNii@##1\DOTS@{\ifx\space@\next\global\DOTS@true\fi}%
  \DN@{\FN@\nextii@}%
  \ifx 3#1\global\DOTSCASE@\z@\else
  \ifx 4#1\global\DOTSCASE@\@ne\else
  \ifx 5#1\global\DOTSCASE@\tw@\else\DN@##1\DOTS@{}%
  \fi\fi\fi\next@}}}
\newif\ifnot@
{\uccode`5=`\\ \uccode`6=`n \uccode`7=`o \uccode`8=`t
 \uppercase{\gdef\not@#1#2#3#4{\relaxnext@
  \DNii@##1\not@{\ifx\space@\next\global\not@true\fi}%
 \global\not@false\DN@##1\not@{}%
 \ifx 5#1\ifx 6#2\ifx 7#3\ifx 8#4\DN@{\FN@\nextii@}\fi\fi\fi
 \fi\next@}}}
\newif\ifkeybin@
\def\keybin@{\keybin@true
 \ifx\next+\else\ifx\next=\else\ifx\next<\else\ifx\next>\else\ifx\next-\else
 \ifx\next*\else\ifx\next:\else\keybin@false\fi\fi\fi\fi\fi\fi\fi}
\def\dots{\RIfM@\expandafter\mdots@\else\expandafter\tdots@\fi}
\def\tdots@{\unskip\relaxnext@
 \DN@{$\m@th\mathinner{\ldotp\ldotp\ldotp}\,
   \ifx\next,\,$\else\ifx\next.\,$\else\ifx\next;\,$\else\ifx\next:\,$\else
   \ifx\next?\,$\else\ifx\next!\,$\else$ \fi\fi\fi\fi\fi\fi}%
 \ \FN@\next@}
\def\mdots@{\FN@\mdots@@}
\def\mdots@@{\gdef\thedots@{\dotso@}
 \ifx\next\boldkey\gdef\thedots@\boldkey{\boldkeydots@}\else                
 \ifx\next\boldsymbol\gdef\thedots@\boldsymbol{\boldsymboldots@}\else       
 \ifx,\next\gdef\thedots@{\dotsc}
 \else\ifx\not\next\gdef\thedots@{\dotsb@}
 \else\keybin@
 \ifkeybin@\gdef\thedots@{\dotsb@}
 \else\xdef\meaning@{\meaning\next..........}\xdef\meaning@@{\meaning@}
  \expandafter\math@\meaning@\math@
  \ifmath@
   \expandafter\mathch@\meaning@\mathch@
   \ifmathch@\expandafter\getmathch@\meaning@\getmathch@\fi                 
  \else\expandafter\macro@\meaning@@\macro@                                 
  \ifmacro@                                                                
   \expandafter\not@\meaning@\not@\ifnot@\gdef\thedots@{\dotsb@}
  \else\expandafter\DOTS@\meaning@\DOTS@
  \ifDOTS@
   \ifcase\number\DOTSCASE@\gdef\thedots@{\dotsb@}%
    \or\gdef\thedots@{\dotsi}\else\fi                                      
  \else\expandafter\math@\meaning@\math@                                   
  \ifmath@\expandafter\mathbin@\meaning@\mathbin@
  \ifmathbin@\gdef\thedots@{\dotsb@}
  \else\expandafter\mathrel@\meaning@\mathrel@
  \ifmathrel@\gdef\thedots@{\dotsb@}
  \fi\fi\fi\fi\fi\fi\fi\fi\fi\fi\fi\fi
 \thedots@}
\def\plainldots@{\mathinner{\ldotp\ldotp\ldotp}}
\def\plaincdots@{\mathinner{\cdotp\cdotp\cdotp}}
\def\dotsi{\!\plaincdots@}
\let\dotsb@\plaincdots@
\newif\ifextra@
\newif\ifrightdelim@
\def\rightdelim@{\global\rightdelim@true                                    
 \ifx\next)\else                                                            
 \ifx\next]\else
 \ifx\next\rbrack\else
 \ifx\next\}\else
 \ifx\next\rbrace\else
 \ifx\next\rangle\else
 \ifx\next\rceil\else
 \ifx\next\rfloor\else
 \ifx\next\rgroup\else
 \ifx\next\rmoustache\else
 \ifx\next\right\else
 \ifx\next\bigr\else
 \ifx\next\biggr\else
 \ifx\next\Bigr\else                                                        
 \ifx\next\Biggr\else\global\rightdelim@false
 \fi\fi\fi\fi\fi\fi\fi\fi\fi\fi\fi\fi\fi\fi\fi}
\def\extra@{%
 \global\extra@false\rightdelim@\ifrightdelim@\global\extra@true            
 \else\ifx\next$\global\extra@true                                          
 \else\xdef\meaning@{\meaning\next..........}
 \expandafter\macro@\meaning@\macro@\ifmacro@                               
 \expandafter\DOTS@\meaning@\DOTS@
 \ifDOTS@
 \ifnum\DOTSCASE@=\tw@\global\extra@true                                    
 \fi\fi\fi\fi\fi}
\newif\ifbold@
\def\dotso@{\relaxnext@
 \ifbold@
  \let\next\delayed@
  \DNii@{\extra@\plainldots@\ifextra@\,\fi}%
 \else
  \DNii@{\DN@{\extra@\plainldots@\ifextra@\,\fi}\FN@\next@}%
 \fi
 \nextii@}
\def\extrap@#1{%
 \ifx\next,\DN@{#1\,}\else
 \ifx\next;\DN@{#1\,}\else
 \ifx\next.\DN@{#1\,}\else\extra@
 \ifextra@\DN@{#1\,}\else
 \let\next@#1\fi\fi\fi\fi\next@}
\def\ldots{\DN@{\extrap@\plainldots@}%
 \FN@\next@}
\def\cdots{\DN@{\extrap@\plaincdots@}%
 \FN@\next@}

\def\dotsc{\relaxnext@
 \DN@{\ifx\next;\plainldots@\,\else
  \ifx\next.\plainldots@\,\else\extra@\plainldots@
  \ifextra@\,\fi\fi\fi}%
 \FN@\next@}
\def\cdot{\mathchar"2201 }

\message{special superscripts,}
\def\dddot#1{{\mathop{#1}\limits^{\vbox to-1.4\ex@{\kern-\tw@\ex@
 \hbox{\rm...}\vss}}}}
\def\ddddot#1{{\mathop{#1}\limits^{\vbox to-1.4\ex@{\kern-\tw@\ex@
 \hbox{\rm....}\vss}}}}
\def\sphat{^{\mathchoice{}{}%
 {\,\,\botsmash{\hbox{\lower4\ex@\hbox{$\m@th\widehat{\null}$}}}}%
 {\,\botsmash{\hbox{\lower3\ex@\hbox{$\m@th\hat{\null}$}}}}}}

\def\spacute{^{\!\botsmash{\hbox{\lower\@ne ex\hbox{\'{}}}}}}
\def\spgrave{^{\mathchoice{}{}{}{\!}%
 \botsmash{\hbox{\lower\@ne ex\hbox{\`{}}}}}}
\def\spdot{^{\hbox{\raise\ex@\hbox{\rm.}}}}
\def\spddot{^{\hbox{\raise\ex@\hbox{\rm..}}}}
\def\spdddot{^{\hbox{\raise\ex@\hbox{\rm...}}}}
\def\spddddot{^{\hbox{\raise\ex@\hbox{\rm....}}}}
\def\spbreve{^{\!\botsmash{\hbox{\lower4\ex@\hbox{\u{}}}}}}

\message{\string\text,}
\def\textonlyfont@#1#2{\def#1{\RIfM@
 \Err@{Use \string#1\space only in text}\else#2\fi}}
\textonlyfont@\rm\tenrm
\textonlyfont@\it\tenit
\textonlyfont@\sl\tensl
\textonlyfont@\bf\tenbf
\def\oldnos#1{\RIfM@{\mathcode`\,="013B \fam\@ne#1}\else
 \leavevmode\hbox{$\m@th\mathcode`\,="013B \fam\@ne#1$}\fi}
\def\text{\RIfM@\expandafter\text@\else\expandafter\text@@\fi}
\def\text@@#1{\leavevmode\hbox{#1}}
\def\mathhexbox@#1#2#3{\text{$\m@th\mathchar"#1#2#3$}}
\def\dag{{\mathhexbox@279}}
\def\ddag{{\mathhexbox@27A}}
\def\S{{\mathhexbox@278}}
\def\P{{\mathhexbox@27B}}
\newif\iffirstchoice@
\firstchoice@true
\def\text@#1{\mathchoice
 {\hbox{\everymath{\displaystyle}\def\textfonti{\the\textfont\@ne}%
  \def\textfontii{\the\textfont\tw@}\textdef@@ T#1}}
 {\hbox{\firstchoice@false
  \everymath{\textstyle}\def\textfonti{\the\textfont\@ne}%
  \def\textfontii{\the\textfont\tw@}\textdef@@ T#1}}
 {\hbox{\firstchoice@false
  \everymath{\scriptstyle}\def\textfonti{\the\scriptfont\@ne}%
  \def\textfontii{\the\scriptfont\tw@}\textdef@@ S\rm#1}}
 {\hbox{\firstchoice@false
  \everymath{\scriptscriptstyle}\def\textfonti
  {\the\scriptscriptfont\@ne}%
  \def\textfontii{\the\scriptscriptfont\tw@}\textdef@@ s\rm#1}}}
\def\textdef@@#1{\textdef@#1\rm\textdef@#1\bf\textdef@#1\sl\textdef@#1\it}
\def\rmfam{0}
\def\textdef@#1#2{%
 \DN@{\csname\expandafter\eat@\string#2fam\endcsname}%
 \if S#1\edef#2{\the\scriptfont\next@\relax}%
 \else\if s#1\edef#2{\the\scriptscriptfont\next@\relax}%
 \else\edef#2{\the\textfont\next@\relax}\fi\fi}
\scriptfont\itfam\tenit \scriptscriptfont\itfam\tenit
\scriptfont\slfam\tensl \scriptscriptfont\slfam\tensl
\newif\iftopfolded@
\newif\ifbotfolded@
\def\topfoldedtext{\topfolded@true\botfolded@false\foldedtext@}
\def\botfoldedtext{\botfolded@true\topfolded@false\foldedtext@}
\def\foldedtext{\topfolded@false\botfolded@false\foldedtext@}
\Invalid@\foldedwidth
\def\foldedtext@{\relaxnext@
 \DN@{\ifx\next\foldedwidth\let\next@\nextii@\else
  \DN@{\nextii@\foldedwidth{.3\hsize}}\fi\next@}%
 \DNii@\foldedwidth##1##2{\setbox\z@\vbox
  {\normalbaselines\hsize##1\relax
  \tolerance1600 \noindent\ignorespaces##2}\ifbotfolded@\boxz@\else
  \iftopfolded@\vtop{\unvbox\z@}\else\vcenter{\boxz@}\fi\fi}%
 \FN@\next@}
\message{math font commands,}
\def\bold{\RIfM@\expandafter\bold@\else
 \expandafter\nonmatherr@\expandafter\bold\fi}
\def\bold@#1{{\bold@@{#1}}}
\def\bold@@#1{\fam\bffam\relax#1}
\def\slanted{\RIfM@\expandafter\slanted@\else
 \expandafter\nonmatherr@\expandafter\slanted\fi}
\def\slanted@#1{{\slanted@@{#1}}}
\def\slanted@@#1{\fam\slfam\relax#1}
\def\roman{\RIfM@\expandafter\roman@\else
 \expandafter\nonmatherr@\expandafter\roman\fi}
\def\roman@#1{{\roman@@{#1}}}
\def\roman@@#1{\fam\rmfam\relax#1}
\def\italic{\RIfM@\expandafter\italic@\else
 \expandafter\nonmatherr@\expandafter\italic\fi}
\def\italic@#1{{\italic@@{#1}}}
\def\italic@@#1{\fam\itfam\relax#1}
\def\Cal{\RIfM@\expandafter\Cal@\else
 \expandafter\nonmatherr@\expandafter\Cal\fi}
\def\Cal@#1{{\Cal@@{#1}}}
\def\Cal@@#1{\noaccents@\fam\tw@#1}
\mathchardef\Gamma="0000
\mathchardef\Delta="0001
\mathchardef\Theta="0002
\mathchardef\Lambda="0003
\mathchardef\Xi="0004
\mathchardef\Pi="0005
\mathchardef\Sigma="0006
\mathchardef\Upsilon="0007
\mathchardef\Phi="0008
\mathchardef\Psi="0009
\mathchardef\Omega="000A
\mathchardef\varGamma="0100
\mathchardef\varDelta="0101
\mathchardef\varTheta="0102
\mathchardef\varLambda="0103
\mathchardef\varXi="0104
\mathchardef\varPi="0105
\mathchardef\varSigma="0106
\mathchardef\varUpsilon="0107
\mathchardef\varPhi="0108
\mathchardef\varPsi="0109
\mathchardef\varOmega="010A
\let\alloc@@\alloc@
\def\hexnumber@#1{\ifcase#1 0\or 1\or 2\or 3\or 4\or 5\or 6\or 7\or 8\or
 9\or A\or B\or C\or D\or E\or F\fi}
\def\loadmsam{%
 \font@\tenmsa=msam10
 \font@\sevenmsa=msam7
 \font@\fivemsa=msam5
 \alloc@@8\fam\chardef\sixt@@n\msafam
 \textfont\msafam=\tenmsa
 \scriptfont\msafam=\sevenmsa
 \scriptscriptfont\msafam=\fivemsa
 \edef\next{\hexnumber@\msafam}%
 \mathchardef\dabar@"0\next39
 \edef\dashrightarrow{\mathrel{\dabar@\dabar@\mathchar"0\next4B}}%
 \edef\dashleftarrow{\mathrel{\mathchar"0\next4C\dabar@\dabar@}}%
 \let\dasharrow\dashrightarrow
 \edef\ulcorner{\delimiter"4\next70\next70 }%
 \edef\urcorner{\delimiter"5\next71\next71 }%
 \edef\llcorner{\delimiter"4\next78\next78 }%
 \edef\lrcorner{\delimiter"5\next79\next79 }%
 \edef\yen{{\noexpand\mathhexbox@\next55}}%
 \edef\checkmark{{\noexpand\mathhexbox@\next58}}%
 \edef\circledR{{\noexpand\mathhexbox@\next72}}%
 \edef\maltese{{\noexpand\mathhexbox@\next7A}}%
 \global\let\loadmsam\empty}%
\def\loadmsbm{%
 \font@\tenmsb=msbm10 \font@\sevenmsb=msbm7 \font@\fivemsb=msbm5
 \alloc@@8\fam\chardef\sixt@@n\msbfam
 \textfont\msbfam=\tenmsb
 \scriptfont\msbfam=\sevenmsb \scriptscriptfont\msbfam=\fivemsb
 \global\let\loadmsbm\empty
 }
\def\widehat#1{\ifx\undefined\msbfam \DN@{362}%
  \else \setboxz@h{$\m@th#1$}%
    \edef\next@{\ifdim\wdz@>\tw@ em%
        \hexnumber@\msbfam 5B%
      \else 362\fi}\fi
  \mathaccent"0\next@{#1}}
\def\widetilde#1{\ifx\undefined\msbfam \DN@{365}%
  \else \setboxz@h{$\m@th#1$}%
    \edef\next@{\ifdim\wdz@>\tw@ em%
        \hexnumber@\msbfam 5D%
      \else 365\fi}\fi
  \mathaccent"0\next@{#1}}
\message{\string\newsymbol,}
\def\newsymbol#1#2#3#4#5{\define#1{}%
  \count@#2\relax \advance\count@\m@ne 
 \ifcase\count@
   \ifx\undefined\msafam\loadmsam\fi \let\next@\msafam
 \or \ifx\undefined\msbfam\loadmsbm\fi \let\next@\msbfam
 \else  \Err@{\Invalid@@\string\newsymbol}\let\next@\tw@\fi
 \mathchardef#1="#3\hexnumber@\next@#4#5\space}

\def\Bbb{\RIfM@\expandafter\Bbb@\else
 \expandafter\nonmatherr@\expandafter\Bbb\fi}
\def\Bbb@#1{{\Bbb@@{#1}}}
\def\Bbb@@#1{\noaccents@\fam\msbfam\relax#1}
\message{bold Greek and bold symbols,}
\def\loadbold{%
 \font@\tencmmib=cmmib10 \font@\sevencmmib=cmmib7 \font@\fivecmmib=cmmib5
 \skewchar\tencmmib'177 \skewchar\sevencmmib'177 \skewchar\fivecmmib'177
 \alloc@@8\fam\chardef\sixt@@n\cmmibfam
 \textfont\cmmibfam\tencmmib
 \scriptfont\cmmibfam\sevencmmib \scriptscriptfont\cmmibfam\fivecmmib
 \font@\tencmbsy=cmbsy10 \font@\sevencmbsy=cmbsy7 \font@\fivecmbsy=cmbsy5
 \skewchar\tencmbsy'60 \skewchar\sevencmbsy'60 \skewchar\fivecmbsy'60
 \alloc@@8\fam\chardef\sixt@@n\cmbsyfam
 \textfont\cmbsyfam\tencmbsy
 \scriptfont\cmbsyfam\sevencmbsy \scriptscriptfont\cmbsyfam\fivecmbsy
 \let\loadbold\empty
}
\def\boldnotloaded#1{\Err@{\ifcase#1\or First\else Second\fi
       bold symbol font not loaded}}
\def\mathchari@#1#2#3{\ifx\undefined\cmmibfam
    \boldnotloaded@\@ne
  \else\mathchar"#1\hexnumber@\cmmibfam#2#3\space \fi}
\def\mathcharii@#1#2#3{\ifx\undefined\cmbsyfam
    \boldnotloaded\tw@
  \else \mathchar"#1\hexnumber@\cmbsyfam#2#3\space\fi}
\edef\bffam@{\hexnumber@\bffam}
\def\boldkey#1{\ifcat\noexpand#1A%
  \ifx\undefined\cmmibfam \boldnotloaded\@ne
  \else {\fam\cmmibfam#1}\fi
 \else
 \ifx#1!\mathchar"5\bffam@21 \else
 \ifx#1(\mathchar"4\bffam@28 \else\ifx#1)\mathchar"5\bffam@29 \else
 \ifx#1+\mathchar"2\bffam@2B \else\ifx#1:\mathchar"3\bffam@3A \else
 \ifx#1;\mathchar"6\bffam@3B \else\ifx#1=\mathchar"3\bffam@3D \else
 \ifx#1?\mathchar"5\bffam@3F \else\ifx#1[\mathchar"4\bffam@5B \else
 \ifx#1]\mathchar"5\bffam@5D \else
 \ifx#1,\mathchari@63B \else
 \ifx#1-\mathcharii@200 \else
 \ifx#1.\mathchari@03A \else
 \ifx#1/\mathchari@03D \else
 \ifx#1<\mathchari@33C \else
 \ifx#1>\mathchari@33E \else
 \ifx#1*\mathcharii@203 \else
 \ifx#1|\mathcharii@06A \else
 \ifx#10\bold0\else\ifx#11\bold1\else\ifx#12\bold2\else\ifx#13\bold3\else
 \ifx#14\bold4\else\ifx#15\bold5\else\ifx#16\bold6\else\ifx#17\bold7\else
 \ifx#18\bold8\else\ifx#19\bold9\else
  \Err@{\string\boldkey\space can't be used with #1}%
 \fi\fi\fi\fi\fi\fi\fi\fi\fi\fi\fi\fi\fi\fi\fi
 \fi\fi\fi\fi\fi\fi\fi\fi\fi\fi\fi\fi\fi\fi}
\def\boldsymbol#1{%
 \DN@{\Err@{You can't use \string\boldsymbol\space with \string#1}#1}%
 \ifcat\noexpand#1A%
   \let\next@\relax
   \ifx\undefined\cmmibfam \boldnotloaded\@ne
   \else {\fam\cmmibfam#1}\fi
 \else
  \xdef\meaning@{\meaning#1.........}%
  \expandafter\math@\meaning@\math@
  \ifmath@
   \expandafter\mathch@\meaning@\mathch@
   \ifmathch@
    \expandafter\boldsymbol@@\meaning@\boldsymbol@@
   \fi
  \else
   \expandafter\macro@\meaning@\macro@
   \expandafter\delim@\meaning@\delim@
   \ifdelim@
    \expandafter\delim@@\meaning@\delim@@
   \else
    \boldsymbol@{#1}%
   \fi
  \fi
 \fi
 \next@}
\def\mathhexboxii@#1#2{\ifx\undefined\cmbsyfam
    \boldnotloaded\tw@
  \else \mathhexbox@{\hexnumber@\cmbsyfam}{#1}{#2}\fi}
\def\boldsymbol@#1{\let\next@\relax\let\next#1%
 \ifx\next\cdot\mathcharii@201 \else
 \ifx\next\prime{{\null\mathcharii@030 \null}}\else
 \ifx\next\lbrack\mathchar"4\bffam@5B \else
 \ifx\next\rbrack\mathchar"5\bffam@5D \else
 \ifx\next\{\mathcharii@466 \else
 \ifx\next\lbrace\mathcharii@466 \else
 \ifx\next\}\mathcharii@567 \else
 \ifx\next\rbrace\mathcharii@567 \else
 \ifx\next\surd{{\mathcharii@170}}\else
 \ifx\next\S{{\mathhexboxii@78}}\else
 \ifx\next\P{{\mathhexboxii@7B}}\else
 \ifx\next\dag{{\mathhexboxii@79}}\else
 \ifx\next\ddag{{\mathhexboxii@7A}}\else
 \DN@{\Err@{You can't use \string\boldsymbol\space with \string#1}#1}%
 \fi\fi\fi\fi\fi\fi\fi\fi\fi\fi\fi\fi\fi}
\def\boldsymbol@@#1.#2\boldsymbol@@{\classnum@#1 \count@@@\classnum@        
 \divide\classnum@4096 \count@\classnum@                                    
 \multiply\count@4096 \advance\count@@@-\count@ \count@@\count@@@           
 \divide\count@@@\@cclvi \count@\count@@                                    
 \multiply\count@@@\@cclvi \advance\count@@-\count@@@                       
 \divide\count@@@\@cclvi                                                    
 \multiply\classnum@4096 \advance\classnum@\count@@                         
 \ifnum\count@@@=\z@                                                        
  \count@"\bffam@ \multiply\count@\@cclvi
  \advance\classnum@\count@
  \DN@{\mathchar\number\classnum@}%
 \else
  \ifnum\count@@@=\@ne                                                      
   \ifx\undefined\cmmibfam \DN@{\boldnotloaded\@ne}%
   \else \count@\cmmibfam \multiply\count@\@cclvi
     \advance\classnum@\count@
     \DN@{\mathchar\number\classnum@}\fi
  \else
   \ifnum\count@@@=\tw@                                                    
     \ifx\undefined\cmbsyfam
       \DN@{\boldnotloaded\tw@}%
     \else
       \count@\cmbsyfam \multiply\count@\@cclvi
       \advance\classnum@\count@
       \DN@{\mathchar\number\classnum@}%
     \fi
  \fi
 \fi
\fi}
\newif\ifdelim@
\newcount\delimcount@
{\uccode`6=`\\ \uccode`7=`d \uccode`8=`e \uccode`9=`l
 \uppercase{\gdef\delim@#1#2#3#4#5\delim@
  {\delim@false\ifx 6#1\ifx 7#2\ifx 8#3\ifx 9#4\delim@true
   \xdef\meaning@{#5}\fi\fi\fi\fi}}}
\def\delim@@#1"#2#3#4#5#6\delim@@{\if#32%
\let\next@\relax
 \ifx\undefined\cmbsyfam \boldnotloaded\@ne
 \else \mathcharii@#2#4#5\space \fi\fi}
\def\vert{\delimiter"026A30C }
\def\Vert{\delimiter"026B30D }
\let\|\Vert
\def\backslash{\delimiter"026E30F }
\def\boldkeydots@#1{\bold@true\let\next=#1\let\delayed@=#1\mdots@@
 \boldkey#1\bold@false}  
\def\boldsymboldots@#1{\bold@true\let\next#1\let\delayed@#1\mdots@@
 \boldsymbol#1\bold@false}
\message{Euler fonts,}

\def\frak{\mathfont@\frak}

\def\loadmathfont#1{%
   \expandafter\font@\csname ten#1\endcsname=#110
   \expandafter\font@\csname seven#1\endcsname=#17
   \expandafter\font@\csname five#1\endcsname=#15
   \edef\next{\noexpand\alloc@@8\fam\chardef\sixt@@n
     \expandafter\noexpand\csname#1fam\endcsname}%
   \next
   \textfont\csname#1fam\endcsname \csname ten#1\endcsname
   \scriptfont\csname#1fam\endcsname \csname seven#1\endcsname
   \scriptscriptfont\csname#1fam\endcsname \csname five#1\endcsname
   \expandafter\def\csname #1\expandafter\endcsname\expandafter{%
      \expandafter\mathfont@\csname#1\endcsname}%
 \expandafter\gdef\csname load#1\endcsname{}%
}
\def\mathfont@#1{\RIfM@\expandafter\mathfont@@\expandafter#1\else
  \expandafter\nonmatherr@\expandafter#1\fi}
\def\mathfont@@#1#2{{\mathfont@@@#1{#2}}}
\def\mathfont@@@#1#2{\noaccents@
   \fam\csname\expandafter\eat@\string#1fam\endcsname
   \relax#2}
\message{math accents,}
\def\accentclass@{7}
\def\noaccents@{\def\accentclass@{0}}
\def\makeacc@#1#2{\def#1{\mathaccent"\accentclass@#2 }}
\makeacc@\hat{05E}
\makeacc@\check{014}
\makeacc@\tilde{07E}
\makeacc@\acute{013}
\makeacc@\grave{012}
\makeacc@\dot{05F}
\makeacc@\ddot{07F}
\makeacc@\breve{015}
\makeacc@\bar{016}

\newcount\skewcharcount@
\newcount\familycount@
\def\theskewchar@{\familycount@\@ne
 \global\skewcharcount@\the\skewchar\textfont\@ne                           
 \ifnum\fam>\m@ne\ifnum\fam<16
  \global\familycount@\the\fam\relax
  \global\skewcharcount@\the\skewchar\textfont\the\fam\relax\fi\fi          
 \ifnum\skewcharcount@>\m@ne
  \ifnum\skewcharcount@<128
  \multiply\familycount@256
  \global\advance\skewcharcount@\familycount@
  \global\advance\skewcharcount@28672
  \mathchar\skewcharcount@\else
  \global\skewcharcount@\m@ne\fi\else
 \global\skewcharcount@\m@ne\fi}                                            
\newcount\pointcount@
\def\getpoints@#1.#2\getpoints@{\pointcount@#1 }
\newdimen\accentdimen@
\newcount\accentmu@
\def\dimentomu@{\multiply\accentdimen@ 100
 \expandafter\getpoints@\the\accentdimen@\getpoints@
 \multiply\pointcount@18
 \divide\pointcount@\@m
 \global\accentmu@\pointcount@}
\def\Makeacc@#1#2{\def#1{\RIfM@\DN@{\mathaccent@
 {"\accentclass@#2 }}\else\DN@{\nonmatherr@{#1}}\fi\next@}}
\def\unbracefonts@{\let\Cal@\Cal@@\let\roman@\roman@@\let\bold@\bold@@
 \let\slanted@\slanted@@}
\def\mathaccent@#1#2{\ifnum\fam=\m@ne\xdef\thefam@{1}\else
 \xdef\thefam@{\the\fam}\fi                                                 
 \accentdimen@\z@                                                           
 \setboxz@h{\unbracefonts@$\m@th\fam\thefam@\relax#2$}
 \ifdim\accentdimen@=\z@\DN@{\mathaccent#1{#2}}
  \setbox@ne\hbox{\unbracefonts@$\m@th\fam\thefam@\relax#2\theskewchar@$}
  \setbox\tw@\hbox{$\m@th\ifnum\skewcharcount@=\m@ne\else
   \mathchar\skewcharcount@\fi$}
  \global\accentdimen@\wd@ne\global\advance\accentdimen@-\wdz@
  \global\advance\accentdimen@-\wd\tw@                                     
  \global\multiply\accentdimen@\tw@
  \dimentomu@\global\advance\accentmu@\@ne                                 
 \else\DN@{{\mathaccent#1{#2\mkern\accentmu@ mu}%
    \mkern-\accentmu@ mu}{}}\fi                                             
 \next@}\Makeacc@\Hat{05E}
\Makeacc@\Check{014}
\Makeacc@\Tilde{07E}
\Makeacc@\Acute{013}
\Makeacc@\Grave{012}
\Makeacc@\Dot{05F}
\Makeacc@\Ddot{07F}
\Makeacc@\Breve{015}
\Makeacc@\Bar{016}
\def\Vec{\RIfM@\DN@{\mathaccent@{"017E }}\else
 \DN@{\nonmatherr@\Vec}\fi\next@}
\def\accentedsymbol#1#2{\csname newbox\expandafter\endcsname
  \csname\expandafter\eat@\string#1@box\endcsname
 \expandafter\setbox\csname\expandafter\eat@
  \string#1@box\endcsname\hbox{$\m@th#2$}\define
  #1{\copy\csname\expandafter\eat@\string#1@box\endcsname{}}}
\message{roots,}
\def\sqrt#1{\radical"270370 {#1}}
\let\underline@\underline
\let\overline@\overline
\def\underline#1{\underline@{#1}}
\def\overline#1{\overline@{#1}}
\Invalid@\leftroot
\Invalid@\uproot
\newcount\uproot@
\newcount\leftroot@
\def\root{\relaxnext@
  \DN@{\ifx\next\uproot\let\next@\nextii@\else
   \ifx\next\leftroot\let\next@\nextiii@\else
   \let\next@\plainroot@\fi\fi\next@}%
  \DNii@\uproot##1{\uproot@##1\relax\FN@\nextiv@}%
  \def\nextiv@{\ifx\next\space@\DN@. {\FN@\nextv@}\else
   \DN@.{\FN@\nextv@}\fi\next@.}%
  \def\nextv@{\ifx\next\leftroot\let\next@\nextvi@\else
   \let\next@\plainroot@\fi\next@}%
  \def\nextvi@\leftroot##1{\leftroot@##1\relax\plainroot@}%
   \def\nextiii@\leftroot##1{\leftroot@##1\relax\FN@\nextvii@}%
  \def\nextvii@{\ifx\next\space@
   \DN@. {\FN@\nextviii@}\else
   \DN@.{\FN@\nextviii@}\fi\next@.}%
  \def\nextviii@{\ifx\next\uproot\let\next@\nextix@\else
   \let\next@\plainroot@\fi\next@}%
  \def\nextix@\uproot##1{\uproot@##1\relax\plainroot@}%
  \bgroup\uproot@\z@\leftroot@\z@\FN@\next@}
\def\plainroot@#1\of#2{\setbox\rootbox\hbox{$\m@th\scriptscriptstyle{#1}$}%
 \mathchoice{\r@@t\displaystyle{#2}}{\r@@t\textstyle{#2}}
 {\r@@t\scriptstyle{#2}}{\r@@t\scriptscriptstyle{#2}}\egroup}
\def\r@@t#1#2{\setboxz@h{$\m@th#1\sqrt{#2}$}%
 \dimen@\ht\z@\advance\dimen@-\dp\z@
 \setbox@ne\hbox{$\m@th#1\mskip\uproot@ mu$}\advance\dimen@ 1.667\wd@ne
 \mkern-\leftroot@ mu\mkern5mu\raise.6\dimen@\copy\rootbox
 \mkern-10mu\mkern\leftroot@ mu\boxz@}
\def\boxed#1{\setboxz@h{$\m@th\displaystyle{#1}$}\dimen@.4\ex@
 \advance\dimen@3\ex@\advance\dimen@\dp\z@
 \hbox{\lower\dimen@\hbox{%
 \vbox{\hrule height.4\ex@
 \hbox{\vrule width.4\ex@\hskip3\ex@\vbox{\vskip3\ex@\boxz@\vskip3\ex@}%
 \hskip3\ex@\vrule width.4\ex@}\hrule height.4\ex@}%
 }}}
\message{commutative diagrams,}
\let\ampersand@\relax
\newdimen\minaw@
\minaw@11.11128\ex@
\newdimen\minCDaw@
\minCDaw@2.5pc
\def\minCDarrowwidth#1{\RIfMIfI@\onlydmatherr@\minCDarrowwidth
 \else\minCDaw@#1\relax\fi\else\onlydmatherr@\minCDarrowwidth\fi}
\newif\ifCD@
\def\CD{\bgroup\vspace@\relax\let\ampersand@&\iffalse}\fi
 \CD@true\vcenter\bgroup\Let@\tabskip\z@skip\baselineskip20\ex@
 \lineskip3\ex@\lineskiplimit3\ex@\halign\bgroup
 &\hfill$\m@th##$\hfill\crcr}
\def\endCD{\crcr\egroup\egroup\egroup}
\newdimen\bigaw@
\atdef@>#1>#2>{\ampersand@                                                  
 \setboxz@h{$\m@th\ssize\;{#1}\;\;$}
 \setbox@ne\hbox{$\m@th\ssize\;{#2}\;\;$}
 \setbox\tw@\hbox{$\m@th#2$}
 \ifCD@\global\bigaw@\minCDaw@\else\global\bigaw@\minaw@\fi                 
 \ifdim\wdz@>\bigaw@\global\bigaw@\wdz@\fi
 \ifdim\wd@ne>\bigaw@\global\bigaw@\wd@ne\fi                                
 \ifCD@\enskip\fi                                                           
 \ifdim\wd\tw@>\z@
  \mathrel{\mathop{\hbox to\bigaw@{\rightarrowfill@\displaystyle}}%
    \limits^{#1}_{#2}}
 \else\mathrel{\mathop{\hbox to\bigaw@{\rightarrowfill@\displaystyle}}%
    \limits^{#1}}\fi                                                        
 \ifCD@\enskip\fi                                                          
 \ampersand@}                                                              
\atdef@<#1<#2<{\ampersand@\setboxz@h{$\m@th\ssize\;\;{#1}\;$}%
 \setbox@ne\hbox{$\m@th\ssize\;\;{#2}\;$}\setbox\tw@\hbox{$\m@th#2$}%
 \ifCD@\global\bigaw@\minCDaw@\else\global\bigaw@\minaw@\fi
 \ifdim\wdz@>\bigaw@\global\bigaw@\wdz@\fi
 \ifdim\wd@ne>\bigaw@\global\bigaw@\wd@ne\fi
 \ifCD@\enskip\fi
 \ifdim\wd\tw@>\z@
  \mathrel{\mathop{\hbox to\bigaw@{\leftarrowfill@\displaystyle}}%
       \limits^{#1}_{#2}}\else
  \mathrel{\mathop{\hbox to\bigaw@{\leftarrowfill@\displaystyle}}%
       \limits^{#1}}\fi
 \ifCD@\enskip\fi\ampersand@}
\begingroup
 \catcode`\~=\active \lccode`\~=`\@
 \lowercase{%
  \global\atdef@)#1)#2){~>#1>#2>}
  \global\atdef@(#1(#2({~<#1<#2<}}
\endgroup
\atdef@ A#1A#2A{\llap{$\m@th\vcenter{\hbox
 {$\ssize#1$}}$}\Big\uparrow\rlap{$\m@th\vcenter{\hbox{$\ssize#2$}}$}&&}
\atdef@ V#1V#2V{\llap{$\m@th\vcenter{\hbox
 {$\ssize#1$}}$}\Big\downarrow\rlap{$\m@th\vcenter{\hbox{$\ssize#2$}}$}&&}
\atdef@={&\enskip\mathrel
 {\vbox{\hrule width\minCDaw@\vskip3\ex@\hrule width
 \minCDaw@}}\enskip&}
\atdef@|{\Big\Vert&&}
\atdef@\vert{\Big\Vert&&}
\def\pretend#1\haswidth#2{\setboxz@h{$\m@th\scriptstyle{#2}$}\hbox
 to\wdz@{\hfill$\m@th\scriptstyle{#1}$\hfill}}
\message{poor man's bold,}
\def\pmb{\RIfM@\expandafter\mathpalette\expandafter\pmb@\else
 \expandafter\pmb@@\fi}
\def\pmb@@#1{\leavevmode\setboxz@h{#1}%
   \dimen@-\wdz@
   \kern-.5\ex@\copy\z@
   \kern\dimen@\kern.25\ex@\raise.4\ex@\copy\z@
   \kern\dimen@\kern.25\ex@\box\z@
}
\def\binrel@@#1{\ifdim\wd2<\z@\mathbin{#1}\else\ifdim\wd\tw@>\z@
 \mathrel{#1}\else{#1}\fi\fi}
\newdimen\pmbraise@
\def\pmb@#1#2{\setbox\thr@@\hbox{$\m@th#1{#2}$}%
 \setbox4\hbox{$\m@th#1\mkern.5mu$}\pmbraise@\wd4\relax
 \binrel@{#2}%
 \dimen@-\wd\thr@@
   \binrel@@{%
   \mkern-.8mu\copy\thr@@
   \kern\dimen@\mkern.4mu\raise\pmbraise@\copy\thr@@
   \kern\dimen@\mkern.4mu\box\thr@@
}}
\def\documentstyle#1{\W@{}\input #1.sty\relax}
\message{syntax check,}
\font\dummyft@=dummy
\fontdimen1 \dummyft@=\z@
\fontdimen2 \dummyft@=\z@
\fontdimen3 \dummyft@=\z@
\fontdimen4 \dummyft@=\z@
\fontdimen5 \dummyft@=\z@
\fontdimen6 \dummyft@=\z@
\fontdimen7 \dummyft@=\z@
\fontdimen8 \dummyft@=\z@
\fontdimen9 \dummyft@=\z@
\fontdimen10 \dummyft@=\z@
\fontdimen11 \dummyft@=\z@
\fontdimen12 \dummyft@=\z@
\fontdimen13 \dummyft@=\z@
\fontdimen14 \dummyft@=\z@
\fontdimen15 \dummyft@=\z@
\fontdimen16 \dummyft@=\z@
\fontdimen17 \dummyft@=\z@
\fontdimen18 \dummyft@=\z@
\fontdimen19 \dummyft@=\z@
\fontdimen20 \dummyft@=\z@
\fontdimen21 \dummyft@=\z@
\fontdimen22 \dummyft@=\z@
\def\fontlist@{\\{\tenrm}\\{\sevenrm}\\{\fiverm}\\{\teni}\\{\seveni}%
 \\{\fivei}\\{\tensy}\\{\sevensy}\\{\fivesy}\\{\tenex}\\{\tenbf}\\{\sevenbf}%
 \\{\fivebf}\\{\tensl}\\{\tenit}}
\def\font@#1=#2 {\rightappend@#1\to\fontlist@\font#1=#2 }
\def\dodummy@{{\def\\##1{\global\let##1\dummyft@}\fontlist@}}
\def\nopages@{\output{\setbox\z@\box\@cclv \deadcycles\z@}%
 \alloc@5\toks\toksdef\@cclvi\output}
\let\galleys\nopages@
\newif\ifsyntax@
\newcount\countxviii@
\def\syntax{\syntax@true\dodummy@\countxviii@\count18
 \loop\ifnum\countxviii@>\m@ne\textfont\countxviii@=\dummyft@
 \scriptfont\countxviii@=\dummyft@\scriptscriptfont\countxviii@=\dummyft@
 \advance\countxviii@\m@ne\repeat                                           
 \dummyft@\tracinglostchars\z@\nopages@\frenchspacing\hbadness\@M}
\def\first@#1#2\end{#1}
\def\printoptions{\W@{Do you want S(yntax check),
  G(alleys) or P(ages)?}%
 \message{Type S, G or P, followed by <return>: }%
 \begingroup 
 \endlinechar\m@ne 
 \read\m@ne to\ans@
 \edef\ans@{\uppercase{\def\noexpand\ans@{%
   \expandafter\first@\ans@ P\end}}}%
 \expandafter\endgroup\ans@
 \if\ans@ P
 \else \if\ans@ S\syntax
 \else \if\ans@ G\galleys
 \else\message{? Unknown option: \ans@; using the `pages' option.}%
 \fi\fi\fi}
\def\alloc@#1#2#3#4#5{\global\advance\count1#1by\@ne
 \ch@ck#1#4#2\allocationnumber=\count1#1
 \global#3#5=\allocationnumber
 \ifalloc@\wlog{\string#5=\string#2\the\allocationnumber}\fi}
\def\document{\def\alloclist@{}\def\fontlist@{}}
\let\enddocument\bye

\let\proclaim\undefined
\let\footnote\undefined
\let\=\undefined
\let\>\undefined

\catcode`\@=\active
\message{... finished}

\expandafter\ifx\csname mathdefs.tex\endcsname\relax
  \expandafter\gdef\csname mathdefs.tex\endcsname{}
\else \message{Hey!  Apparently you were trying to
  \string\input{mathdefs.tex} twice.   This does not make sense.} 
\errmessage{Please edit your file (probably \jobname.tex) and remove
any duplicate ``\string\input'' lines}\endinput\fi




\catcode`\X=12\catcode`\@=11

\def\n@wcount{\alloc@0\count\countdef\insc@unt}
\def\n@wwrite{\alloc@7\write\chardef\sixt@@n}
\def\n@wread{\alloc@6\read\chardef\sixt@@n}
\def\r@s@t{\relax}\def\v@idline{\par}\def\@mputate#1/{#1}
\def\l@c@l#1X{\firstpart.#1}\def\gl@b@l#1X{#1}\def\t@d@l#1X{{}}

\def\crossrefs#1{\ifx\all#1\let\tr@ce=\all\else\def\tr@ce{#1,}\fi
   \n@wwrite\cit@tionsout\openout\cit@tionsout=\jobname.cit 
   \write\cit@tionsout{\tr@ce}\expandafter\setfl@gs\tr@ce,}
\def\setfl@gs#1,{\def\@{#1}\ifx\@\empty\let\next=\relax
   \else\let\next=\setfl@gs\expandafter\xdef
   \csname#1tr@cetrue\endcsname{}\fi\next}
\def\m@ketag#1#2{\expandafter\n@wcount\csname#2tagno\endcsname
     \csname#2tagno\endcsname=0\let\tail=\all\xdef\all{\tail#2,}
   \ifx#1\l@c@l\let\tail=\r@s@t\xdef\r@s@t{\csname#2tagno\endcsname=0\tail}\fi
   \expandafter\gdef\csname#2cite\endcsname##1{\expandafter
     \ifx\csname#2tag##1\endcsname\relax?\else\csname#2tag##1\endcsname\fi
     \expandafter\ifx\csname#2tr@cetrue\endcsname\relax\else
     \write\cit@tionsout{#2tag ##1 cited on page \folio.}\fi}
   \expandafter\gdef\csname#2page\endcsname##1{\expandafter
     \ifx\csname#2page##1\endcsname\relax?\else\csname#2page##1\endcsname\fi
     \expandafter\ifx\csname#2tr@cetrue\endcsname\relax\else
     \write\cit@tionsout{#2tag ##1 cited on page \folio.}\fi}
   \expandafter\gdef\csname#2tag\endcsname##1{\expandafter
      \ifx\csname#2check##1\endcsname\relax
      \expandafter\xdef\csname#2check##1\endcsname{}%
      \else\immediate\write16{Warning: #2tag ##1 used more than once.}\fi
      \multit@g{#1}{#2}##1/X%
      \write\t@gsout{#2tag ##1 assigned number \csname#2tag##1\endcsname\space
      on page \number\count0.}%
   \csname#2tag##1\endcsname}}

\def\multit@g#1#2#3/#4X{\def\t@mp{#4}\ifx\t@mp\empty%
      \global\advance\csname#2tagno\endcsname by 1 
      \expandafter\xdef\csname#2tag#3\endcsname
      {#1\number\csname#2tagno\endcsnameX}%
   \else\expandafter\ifx\csname#2last#3\endcsname\relax
      \expandafter\n@wcount\csname#2last#3\endcsname
      \global\advance\csname#2tagno\endcsname by 1 
      \expandafter\xdef\csname#2tag#3\endcsname
      {#1\number\csname#2tagno\endcsnameX}
      \write\t@gsout{#2tag #3 assigned number \csname#2tag#3\endcsname\space
      on page \number\count0.}\fi
   \global\advance\csname#2last#3\endcsname by 1
   \def\t@mp{\expandafter\xdef\csname#2tag#3/}%
   \expandafter\t@mp\@mputate#4\endcsname
   {\csname#2tag#3\endcsname\lastpart{\csname#2last#3\endcsname}}\fi}
\def\t@gs#1{\def\all{}\m@ketag#1e\m@ketag#1s\m@ketag\t@d@l p
\let\realscite\scite
\let\realstag\stag
   \m@ketag\gl@b@l r \n@wread\t@gsin
   \openin\t@gsin=\jobname.tgs \re@der \closein\t@gsin
   \n@wwrite\t@gsout\openout\t@gsout=\jobname.tgs }
\outer\def\localtags{\t@gs\l@c@l}
\outer\def\globaltags{\t@gs\gl@b@l}
\outer\def\newlocaltag#1{\m@ketag\l@c@l{#1}}
\outer\def\newglobaltag#1{\m@ketag\gl@b@l{#1}}

\newif\ifpr@ 
\def\m@kecs #1tag #2 assigned number #3 on page #4.%
   {\expandafter\gdef\csname#1tag#2\endcsname{#3}
   \expandafter\gdef\csname#1page#2\endcsname{#4}
   \ifpr@\expandafter\xdef\csname#1check#2\endcsname{}\fi}
\def\re@der{\ifeof\t@gsin\let\next=\relax\else
   \read\t@gsin to\t@gline\ifx\t@gline\v@idline\else
   \expandafter\m@kecs \t@gline\fi\let \next=\re@der\fi\next}
\def\pretags#1{\pr@true\pret@gs#1,,}
\def\pret@gs#1,{\def\@{#1}\ifx\@\empty\let\n@xtfile=\relax
   \else\let\n@xtfile=\pret@gs \openin\t@gsin=#1.tgs \message{#1} \re@der 
   \closein\t@gsin\fi \n@xtfile}

\newcount\sectno\sectno=0\newcount\subsectno\subsectno=0
\newif\ifultr@local \def\ultralocal{\ultr@localtrue}
\def\firstpart{\number\sectno}
\def\lastpart#1{\ifcase#1 \or a\or b\or c\or d\or e\or f\or g\or h\or 
   i\or k\or l\or m\or n\or o\or p\or q\or r\or s\or t\or u\or v\or w\or 
   x\or y\or z \fi}

\def\resetall{\global\advance\sectno by 1\subsectno=0
   \gdef\firstpart{\number\sectno}\r@s@t}
\def\resetsub{\global\advance\subsectno by 1
   \gdef\firstpart{\number\sectno.\number\subsectno}\r@s@t}
\def\newsection#1\par{\resetall\vskip0pt plus.3\vsize\penalty-250
   \vskip0pt plus-.3\vsize\bigskip\bigskip
   \message{#1}\leftline{\bf#1}\nobreak\bigskip}
\def\subsection#1\par{\ifultr@local\resetsub\fi
   \vskip0pt plus.2\vsize\penalty-250\vskip0pt plus-.2\vsize
   \bigskip\smallskip\message{#1}\leftline{\bf#1}\nobreak\medskip}


\newdimen\marginshift

\newdimen\margindelta
\newdimen\marginmax
\newdimen\marginmin

\def\margininit{       
\marginmax=3 true cm                  
				      
\margindelta=0.1 true cm              
\marginmin=0.1true cm                 
\marginshift=\marginmin
}    

\def\t@gsjj#1,{\def\@{#1}\ifx\@\empty\let\next=\relax\else\let\next=\t@gsjj
   \def\@@{p}\ifx\@\@@\else
   \expandafter\gdef\csname#1cite\endcsname##1{\citejj{##1}}
   \expandafter\gdef\csname#1page\endcsname##1{?}
   \expandafter\gdef\csname#1tag\endcsname##1{\tagjj{##1}}\fi\fi\next}
\newif\ifshowstuffinmargin
\showstuffinmarginfalse
\def\jjtags{\ifx\shlhetal\relax 
  \else
\ifx\shlhetal\undefinedcontrolseq
\else
\showstuffinmargintrue
\ifx\all\relax\else\expandafter\t@gsjj\all,\fi\fi \fi
}

\def\tagjj#1{\realstag{#1}\mginpar{\zeigen{#1}}}
\def\citejj#1{\rechnen{#1}\mginpar{\zeigen{#1}}}     

\def\rechnen#1{\expandafter\ifx\csname stag#1\endcsname\relax ??\else
                           \csname stag#1\endcsname\fi}

\newdimen\theight

\def\marginfont{\sevenrm}

\def\trymarginbox#1{\setbox0=\hbox{\marginfont\hskip\marginshift #1}%
		\global\marginshift\wd0 
		\global\advance\marginshift\margindelta}

\def \mginpar#1{%
\ifvmode\setbox0\hbox to \hsize{\hfill\rlap{\marginfont\quad#1}}%
\ht0 0cm
\dp0 0cm
\box0\vskip-\baselineskip
\else 
             \vadjust{\trymarginbox{#1}%
		\ifdim\marginshift>\marginmax \global\marginshift\marginmin
			\trymarginbox{#1}%
                \fi
             \theight=\ht0
             \advance\theight by \dp0    \advance\theight by \lineskip
             \kern -\theight \vbox to \theight{\rightline{\rlap{\box0}}%
\vss}}\fi}


\def\t@gsoff#1,{\def\@{#1}\ifx\@\empty\let\next=\relax\else\let\next=\t@gsoff
   \def\@@{p}\ifx\@\@@\else
   \expandafter\gdef\csname#1cite\endcsname##1{\zeigen{##1}}
   \expandafter\gdef\csname#1page\endcsname##1{?}
   \expandafter\gdef\csname#1tag\endcsname##1{\zeigen{##1}}\fi\fi\next}
\def\verbatimtags{\showstuffinmarginfalse
\ifx\all\relax\else\expandafter\t@gsoff\all,\fi}
\def\zeigen#1{\hbox{$\langle$}#1\hbox{$\rangle$}}

\def\margincite#1{\ifshowstuffinmargin\mginpar{\zeigen{#1}}\fi}

\def\(#1){\edef\dot@g{\ifmmode\ifinner(\hbox{\noexpand\etag{#1}})
   \else\noexpand\eqno(\hbox{\noexpand\etag{#1}})\fi
   \else(\noexpand\ecite{#1})\fi}\dot@g}

\newif\ifbr@ck
\def\eat#1{}
\def\[#1]{\br@cktrue[\br@cket#1'X]}
\def\br@cket#1'#2X{\def\temp{#2}\ifx\temp\empty\let\next\eat
   \else\let\next\br@cket\fi
   \ifbr@ck\br@ckfalse\br@ck@t#1,X\else\br@cktrue#1\fi\next#2X}
\def\br@ck@t#1,#2X{\def\temp{#2}\ifx\temp\empty\let\neext\eat
   \else\let\neext\br@ck@t\def\temp{,}\fi
   \def\teemp{#1}\ifx\teemp\empty\else\rcite{#1}\fi\temp\neext#2X}
\def\resetbr@cket{\gdef\[##1]{[\rtag{##1}]}}
\def\references{\resetbr@cket\newsection References\par}

\newtoks\symb@ls\newtoks\s@mb@ls\newtoks\p@gelist\n@wcount\ftn@mber
    \ftn@mber=1\newif\ifftn@mbers\ftn@mbersfalse\newif\ifbyp@ge\byp@gefalse
\def\defm@rk{\ifftn@mbers\n@mberm@rk\else\symb@lm@rk\fi}
\def\n@mberm@rk{\xdef\m@rk{{\the\ftn@mber}}%
    \global\advance\ftn@mber by 1 }
\def\rot@te#1{\let\temp=#1\global#1=\expandafter\r@t@te\the\temp,X}
\def\r@t@te#1,#2X{{#2#1}\xdef\m@rk{{#1}}}
\def\b@@st#1{{$^{#1}$}}\def\str@p#1{#1}
\def\symb@lm@rk{\ifbyp@ge\rot@te\p@gelist\ifnum\expandafter\str@p\m@rk=1 
    \s@mb@ls=\symb@ls\fi\write\f@nsout{\number\count0}\fi \rot@te\s@mb@ls}
\def\byp@ge{\byp@getrue\n@wwrite\f@nsin\openin\f@nsin=\jobname.fns 
    \n@wcount\currentp@ge\currentp@ge=0\p@gelist={0}
    \re@dfns\closein\f@nsin\rot@te\p@gelist
    \n@wread\f@nsout\openout\f@nsout=\jobname.fns }
\def\m@kelist#1X#2{{#1,#2}}
\def\re@dfns{\ifeof\f@nsin\let\next=\relax\else\read\f@nsin to \f@nline
    \ifx\f@nline\v@idline\else\let\t@mplist=\p@gelist
    \ifnum\currentp@ge=\f@nline
    \global\p@gelist=\expandafter\m@kelist\the\t@mplistX0
    \else\currentp@ge=\f@nline
    \global\p@gelist=\expandafter\m@kelist\the\t@mplistX1\fi\fi
    \let\next=\re@dfns\fi\next}
\def\symbols#1{\symb@ls={#1}\s@mb@ls=\symb@ls} 
\def\bigsymbol{\textstyle}
\symbols{\bigsymbol\ast,\dagger,\ddagger,\sharp,\flat,\natural,\star}
\def\ftnumbers{\ftn@mberstrue} \def\ftsymbols{\ftn@mbersfalse}
\def\paginal{\byp@ge} \def\resetftnumbers{\ftn@mber=1}
\def\ftnote#1{\defm@rk\expandafter\expandafter\expandafter\footnote
    \expandafter\b@@st\m@rk{#1}}

\long\def\jump#1\endjump{}
\def\ssum{\mathop{\lower .1em\hbox{$\textstyle\Sigma$}}\nolimits}

\def\qed{\nobreak\kern 1em \vrule height .5em width .5em depth 0em}
\def\newneq{\hbox{\rlap{\hbox to 1\wd9{\hss$=$\hss}}\raise .1em 
   \hbox to 1\wd9{\hss$\scriptscriptstyle/$\hss}}}
\def\subsetne{\setbox9 = \hbox{$\subset$}\mathrel{\hbox{\rlap
   {\lower .4em \newneq}\raise .13em \hbox{$\subset$}}}}
\def\supsetne{\setbox9 = \hbox{$\subset$}\mathrel{\hbox{\rlap
   {\lower .4em \newneq}\raise .13em \hbox{$\supset$}}}}

\def\vbar{\mathchoice{\vrule height6.3ptdepth-.5ptwidth.8pt\kern-.8pt}
   {\vrule height6.3ptdepth-.5ptwidth.8pt\kern-.8pt}
   {\vrule height4.1ptdepth-.35ptwidth.6pt\kern-.6pt}
   {\vrule height3.1ptdepth-.25ptwidth.5pt\kern-.5pt}}
\def\f@dge{\mathchoice{}{}{\mkern.5mu}{\mkern.8mu}}
\def\b@c#1#2{{\rm \mkern#2mu\vbar\mkern-#2mu#1}}
\def\b@b#1{{\rm I\mkern-3.5mu #1}}
\def\b@a#1#2{{\rm #1\mkern-#2mu\f@dge #1}}
\def\bb#1{{\count4=`#1 \advance\count4by-64 \ifcase\count4\or\b@a A{11.5}\or
   \b@b B\or\b@c C{5}\or\b@b D\or\b@b E\or\b@b F \or\b@c G{5}\or\b@b H\or
   \b@b I\or\b@c J{3}\or\b@b K\or\b@b L \or\b@b M\or\b@b N\or\b@c O{5} \or
   \b@b P\or\b@c Q{5}\or\b@b R\or\b@a S{8}\or\b@a T{10.5}\or\b@c U{5}\or
   \b@a V{12}\or\b@a W{16.5}\or\b@a X{11}\or\b@a Y{11.7}\or\b@a Z{7.5}\fi}}

\catcode`\X=11 \catcode`\@=12




\let\thischap\jobname

\def\partof#1{\csname returnthe#1part\endcsname}
\def\chapof#1{\csname returnthe#1chap\endcsname}

\def\setchapter#1,#2,#3.{%
  \expandafter\def\csname returnthe#1part\endcsname{#2}%
  \expandafter\def\csname returnthe#1chap\endcsname{#3}%
}

\setchapter 300a,A,I.
\setchapter 300b,A,II.
\setchapter 300c,A,III.
\setchapter 300d,A,IV.
\setchapter 300e,A,V.
\setchapter 300f,A,VI.
\setchapter 300g,A,VII.
\setchapter F604,B,0.
\setchapter  88r,B,I.
\setchapter  600,B,II.
\setchapter  705,B,III.

\def\cprefix#1{
\edef\theotherpart{\partof{#1}}\edef\theotherchap{\chapof{#1}}%
\ifx\theotherpart\thispart
   \ifx\theotherchap\thischap 
    \else 
     \theotherchap%
    \fi
   \else 
     \theotherpart.\theotherchap\fi}

\def\sectioncite[#1]#2{%
     \cprefix{#2}#1}

\edef\thispart{\partof{\thischap}}
\edef\thischap{\chapof{\thischap}}


\def\spuriousreset{}


\expandafter\ifx\csname citeadd.tex\endcsname\relax
\expandafter\gdef\csname citeadd.tex\endcsname{}
\else \message{Hey!  Apparently you were trying to
\string\input{citeadd.tex} twice.   This does not make sense.} 
\errmessage{Please edit your file (probably \jobname.tex) and remove
any duplicate ``\string\input'' lines}\endinput\fi

\sectno=-1   
\localtags
\jjtags
\NoBlackBoxes
\define\mr{\medskip\roster}
\define\sn{\smallskip\noindent}
\define\mn{\medskip\noindent}

\define\ub{\underbar}
\define\wilog{\text{without loss of generality}}
\define\ermn{\endroster\medskip\noindent}
\define\dbca{\dsize\bigcap}

\define \nl{\newline}
\magnification=\magstep 1
\documentstyle{amsppt}

{    
\catcode`@11

\ifx\alicetwothousandloaded@\relax
  \endinput\else\global\let\alicetwothousandloaded@\relax\fi

\gdef\subjclass{\let\savedef@\subjclass
 \def\subjclass##1\endsubjclass{\let\subjclass\savedef@
   \toks@{\def\usualspace{{\rm\enspace}}\eightpoint}%
   \toks@@{##1\unskip.}%
   \edef\thesubjclass@{\the\toks@
     \frills@{{\noexpand\rm2000 {\noexpand\it Mathematics Subject
       Classification}.\noexpand\enspace}}%
     \the\toks@@}}%
  \nofrillscheck\subjclass}
} 


\expandafter\ifx\csname alice2jlem.tex\endcsname\relax
  \expandafter\xdef\csname alice2jlem.tex\endcsname{\the\catcode`@}
\else \message{Hey!  Apparently you were trying to
\string\input{alice2jlem.tex}  twice.   This does not make sense.}
\errmessage{Please edit your file (probably \jobname.tex) and remove
any duplicate ``\string\input'' lines}\endinput\fi

\expandafter\ifx\csname bib4plain.tex\endcsname\relax
  \expandafter\gdef\csname bib4plain.tex\endcsname{}
\else \message{Hey!  Apparently you were trying to \string\input
  bib4plain.tex twice.   This does not make sense.}
\errmessage{Please edit your file (probably \jobname.tex) and remove
any duplicate ``\string\input'' lines}\endinput\fi

\def\renewcommand{\newcommand}	       
\edef\cite{\the\catcode`@}%
\catcode`@ = 11
\let\@oldatcatcode = \cite
\chardef\@letter = 11
\chardef\@other = 12
%
%
%
%
\def\@innerdef#1#2{\edef#1{\expandafter\noexpand\csname #2\endcsname}}%
%
%
\@innerdef\@innernewcount{newcount}%
\@innerdef\@innernewdimen{newdimen}%
\@innerdef\@innernewif{newif}%
\@innerdef\@innernewwrite{newwrite}%
%
%
%
\def\@gobble#1{}%
%
%
%
\ifx\inputlineno\@undefined
   \let\@linenumber = \empty 
\else
   \def\@linenumber{\the\inputlineno:\space}%
\fi
%
%
%
\def\@futurenonspacelet#1{\def\cs{#1}%
   \afterassignment\@stepone\let\@nexttoken=
}%
\begingroup 
\def\\{\global\let\@stoken= }%
\\ 
\endgroup
\def\@stepone{\expandafter\futurelet\cs\@steptwo}%
\def\@steptwo{\expandafter\ifx\cs\@stoken\let\@@next=\@stepthree
   \else\let\@@next=\@nexttoken\fi \@@next}%
\def\@stepthree{\afterassignment\@stepone\let\@@next= }%
%
%
%
\def\@getoptionalarg#1{%
   \let\@optionaltemp = #1%
   \let\@optionalnext = \relax
   \@futurenonspacelet\@optionalnext\@bracketcheck
}%
%
%
\def\@bracketcheck{%
   \ifx [\@optionalnext
      \expandafter\@@getoptionalarg
   \else
      \let\@optionalarg = \empty
      \expandafter\@optionaltemp
   \fi
}%
\def\@@getoptionalarg[#1]{%
   \def\@optionalarg{#1}%
   \@optionaltemp
}%
%
%
%
\def\@nnil{\@nil}%
\def\@fornoop#1\@@#2#3{}%
\def\@for#1:=#2\do#3{%
   \edef\@fortmp{#2}%
   \ifx\@fortmp\empty \else
      \expandafter\@forloop#2,\@nil,\@nil\@@#1{#3}%
   \fi
}%
\def\@forloop#1,#2,#3\@@#4#5{\def#4{#1}\ifx #4\@nnil \else
       #5\def#4{#2}\ifx #4\@nnil \else#5\@iforloop #3\@@#4{#5}\fi\fi
}%
\def\@iforloop#1,#2\@@#3#4{\def#3{#1}\ifx #3\@nnil
       \let\@nextwhile=\@fornoop \else
      #4\relax\let\@nextwhile=\@iforloop\fi\@nextwhile#2\@@#3{#4}%
}%
%
%
%
\@innernewif\if@fileexists
\def\@testfileexistence{\@getoptionalarg\@finishtestfileexistence}%
\def\@finishtestfileexistence#1{%
   \begingroup
      \def\extension{#1}%
      \immediate\openin0 =
         \ifx\@optionalarg\empty\jobname\else\@optionalarg\fi
         \ifx\extension\empty \else .#1\fi
         \space
      \ifeof 0
         \global\@fileexistsfalse
      \else
         \global\@fileexiststrue
      \fi
      \immediate\closein0
   \endgroup
}%
%
%
%
%
\def\bibliographystyle#1{%
   \@readauxfile
   \@writeaux{\string\bibstyle{#1}}%
}%
\let\bibstyle = \@gobble
%
%
\let\bblfilebasename = \jobname
\def\bibliography#1{%
   \@readauxfile
   \@writeaux{\string\bibdata{#1}}%
   \@testfileexistence[\bblfilebasename]{bbl}%
   \if@fileexists
      \nobreak
      \@readbblfile
   \fi
}%
\let\bibdata = \@gobble
%
%
\def\nocite#1{%
   \@readauxfile
   \@writeaux{\string\citation{#1}}%
}%
\@innernewif\if@notfirstcitation
%
%
\def\cite{\@getoptionalarg\@cite}%
%
%
\def\@cite#1{%
   \let\@citenotetext = \@optionalarg
   \printcitestart
   \nocite{#1}%
   \@notfirstcitationfalse
   \@for \@citation :=#1\do
   {%
      \expandafter\@onecitation\@citation\@@
   }%
   \ifx\empty\@citenotetext\else
      \printcitenote{\@citenotetext}%
   \fi
   \printcitefinish
}%
\newif\ifweareinprivate
\weareinprivatetrue
\ifx\shlhetal\undefinedcontrolseq\weareinprivatefalse\fi
\ifx\shlhetal\relax\weareinprivatefalse\fi
\def\@onecitation#1\@@{%
   \if@notfirstcitation
      \printbetweencitations
   \fi
   \expandafter \ifx \csname\@citelabel{#1}\endcsname \relax
      \if@citewarning
         \message{\@linenumber Undefined citation `#1'.}%
      \fi
     \ifweareinprivate
      \expandafter\gdef\csname\@citelabel{#1}\endcsname{%
\strut 
\vadjust{\vskip-\dp\strutbox
\vbox to 0pt{\vss\parindent0cm \leftskip=\hsize 
\advance\leftskip3mm
\advance\hsize 4cm\strut\openup-4pt 
\rightskip 0cm plus 1cm minus 0.5cm ?  #1 ?\strut}}
         {\tt
            \escapechar = -1
            \nobreak\hskip0pt\pfeilsw
            \expandafter\string\csname#1\endcsname
             \pfeilso
            \nobreak\hskip0pt
         }%
      }%
     \else  
      \expandafter\gdef\csname\@citelabel{#1}\endcsname{%
            {\tt\expandafter\string\csname#1\endcsname}
      }%
     \fi  
   \fi
   \csname\@citelabel{#1}\endcsname
   \@notfirstcitationtrue
}%
%
%
\def\@citelabel#1{b@#1}%
%
%
\def\@citedef#1#2{\expandafter\gdef\csname\@citelabel{#1}\endcsname{#2}}%
%
%
%
\def\@readbblfile{%
   \ifx\@itemnum\@undefined
      \@innernewcount\@itemnum
   \fi
   \begingroup
      \def\begin##1##2{%
         \setbox0 = \hbox{\biblabelcontents{##2}}%
         \biblabelwidth = \wd0
      }%
      \def\end##1{}
      %
      %
      \@itemnum = 0
      \def\bibitem{\@getoptionalarg\@bibitem}%
      \def\@bibitem{%
         \ifx\@optionalarg\empty
            \expandafter\@numberedbibitem
         \else
            \expandafter\@alphabibitem
         \fi
      }%
      \def\@alphabibitem##1{%
         \expandafter \xdef\csname\@citelabel{##1}\endcsname {\@optionalarg}%
         \ifx\biblabelprecontents\@undefined
            \let\biblabelprecontents = \relax
         \fi
         \ifx\biblabelpostcontents\@undefined
            \let\biblabelpostcontents = \hss
         \fi
         \@finishbibitem{##1}%
      }%
      \def\@numberedbibitem##1{%
         \advance\@itemnum by 1
         \expandafter \xdef\csname\@citelabel{##1}\endcsname{\number\@itemnum}%
         \ifx\biblabelprecontents\@undefined
            \let\biblabelprecontents = \hss
         \fi
         \ifx\biblabelpostcontents\@undefined
            \let\biblabelpostcontents = \relax
         \fi
         \@finishbibitem{##1}%
      }%
      \def\@finishbibitem##1{%
         \biblabelprint{\csname\@citelabel{##1}\endcsname}%
         \@writeaux{\string\@citedef{##1}{\csname\@citelabel{##1}\endcsname}}%
         \ignorespaces
      }%
      %
      %
      \let\em = \bblem
      \let\newblock = \bblnewblock
      \let\sc = \bblsc
      \frenchspacing
      \clubpenalty = 4000 \widowpenalty = 4000
      \tolerance = 10000 \hfuzz = .5pt
      \everypar = {\hangindent = \biblabelwidth
                      \advance\hangindent by \biblabelextraspace}%
      \bblrm
      \parskip = 1.5ex plus .5ex minus .5ex
      \biblabelextraspace = .5em
      \bblhook
      \input \bblfilebasename.bbl
   \endgroup
}%
%
%
\@innernewdimen\biblabelwidth
\@innernewdimen\biblabelextraspace
%
%
%
\def\biblabelprint#1{%
   \noindent
   \hbox to \biblabelwidth{%
      \biblabelprecontents
      \biblabelcontents{#1}%
      \biblabelpostcontents
   }%
   \kern\biblabelextraspace
}%
%
%
%
\def\biblabelcontents#1{{\bblrm [#1]}}%
%
%
\def\bblrm{\rm}%
%
%
\def\bblem{\it}%
%
%
\def\bblsc{\ifx\@scfont\@undefined
              \font\@scfont = cmcsc10
           \fi
           \@scfont
}%
%
%
\def\bblnewblock{\hskip .11em plus .33em minus .07em }%
%
%
\let\bblhook = \empty
%
%
%
\def\printcitestart{[}
\def\printcitefinish{]}
\def\printbetweencitations{, }
\def\printcitenote#1{, #1}
%
%
%
\let\citation = \@gobble
%
%
%
\@innernewcount\@numparams
%
%
\def\newcommand#1{%
   \def\@commandname{#1}%
   \@getoptionalarg\@continuenewcommand
}%
%
%
\def\@continuenewcommand{%
   \@numparams = \ifx\@optionalarg\empty 0\else\@optionalarg \fi \relax
   \@newcommand
}%
%
%
\def\@newcommand#1{%
   \def\@startdef{\expandafter\edef\@commandname}%
   \ifnum\@numparams=0
      \let\@paramdef = \empty
   \else
      \ifnum\@numparams>9
         \errmessage{\the\@numparams\space is too many parameters}%
      \else
         \ifnum\@numparams<0
            \errmessage{\the\@numparams\space is too few parameters}%
         \else
            \edef\@paramdef{%
               \ifcase\@numparams
                  \empty  No arguments.
               \or ####1%
               \or ####1####2%
               \or ####1####2####3%
               \or ####1####2####3####4%
               \or ####1####2####3####4####5%
               \or ####1####2####3####4####5####6%
               \or ####1####2####3####4####5####6####7%
               \or ####1####2####3####4####5####6####7####8%
               \or ####1####2####3####4####5####6####7####8####9%
               \fi
            }%
         \fi
      \fi
   \fi
   \expandafter\@startdef\@paramdef{#1}%
}%
%
%
%
%
\def\@readauxfile{%
   \if@auxfiledone \else 
      \global\@auxfiledonetrue
      \@testfileexistence{aux}%
      \if@fileexists
         \begingroup
            \endlinechar = -1
            \catcode`@ = 11
            \input \jobname.aux
         \endgroup
      \else
         \message{\@undefinedmessage}%
         \global\@citewarningfalse
      \fi
      \immediate\openout\@auxfile = \jobname.aux
   \fi
}%
%
%
\newif\if@auxfiledone
\ifx\noauxfile\@undefined \else \@auxfiledonetrue\fi
%
%
%
%
\@innernewwrite\@auxfile
\def\@writeaux#1{\ifx\noauxfile\@undefined \write\@auxfile{#1}\fi}%
%
%
%
\ifx\@undefinedmessage\@undefined
   \def\@undefinedmessage{No .aux file; I won't give you warnings about
                          undefined citations.}%
\fi
%
%
\@innernewif\if@citewarning
\ifx\noauxfile\@undefined \@citewarningtrue\fi
%
%
%
\catcode`@ = \@oldatcatcode

\def\pfeilso{\leavevmode
            \vrule width 1pt height9pt depth 0pt\relax
           \vrule width 1pt height8.7pt depth 0pt\relax
           \vrule width 1pt height8.3pt depth 0pt\relax
           \vrule width 1pt height8.0pt depth 0pt\relax
           \vrule width 1pt height7.7pt depth 0pt\relax
            \vrule width 1pt height7.3pt depth 0pt\relax
            \vrule width 1pt height7.0pt depth 0pt\relax
            \vrule width 1pt height6.7pt depth 0pt\relax
            \vrule width 1pt height6.3pt depth 0pt\relax
            \vrule width 1pt height6.0pt depth 0pt\relax
            \vrule width 1pt height5.7pt depth 0pt\relax
            \vrule width 1pt height5.3pt depth 0pt\relax
            \vrule width 1pt height5.0pt depth 0pt\relax
            \vrule width 1pt height4.7pt depth 0pt\relax
            \vrule width 1pt height4.3pt depth 0pt\relax
            \vrule width 1pt height4.0pt depth 0pt\relax
            \vrule width 1pt height3.7pt depth 0pt\relax
            \vrule width 1pt height3.3pt depth 0pt\relax
            \vrule width 1pt height3.0pt depth 0pt\relax
            \vrule width 1pt height2.7pt depth 0pt\relax
            \vrule width 1pt height2.3pt depth 0pt\relax
            \vrule width 1pt height2.0pt depth 0pt\relax
            \vrule width 1pt height1.7pt depth 0pt\relax
            \vrule width 1pt height1.3pt depth 0pt\relax
            \vrule width 1pt height1.0pt depth 0pt\relax
            \vrule width 1pt height0.7pt depth 0pt\relax
            \vrule width 1pt height0.3pt depth 0pt\relax}

\def\pfeilsw{ \leavevmode 
            \vrule width 1pt height0.3pt depth 0pt\relax
            \vrule width 1pt height0.7pt depth 0pt\relax
            \vrule width 1pt height1.0pt depth 0pt\relax
            \vrule width 1pt height1.3pt depth 0pt\relax
            \vrule width 1pt height1.7pt depth 0pt\relax
            \vrule width 1pt height2.0pt depth 0pt\relax
            \vrule width 1pt height2.3pt depth 0pt\relax
            \vrule width 1pt height2.7pt depth 0pt\relax
            \vrule width 1pt height3.0pt depth 0pt\relax
            \vrule width 1pt height3.3pt depth 0pt\relax
            \vrule width 1pt height3.7pt depth 0pt\relax
            \vrule width 1pt height4.0pt depth 0pt\relax
            \vrule width 1pt height4.3pt depth 0pt\relax
            \vrule width 1pt height4.7pt depth 0pt\relax
            \vrule width 1pt height5.0pt depth 0pt\relax
            \vrule width 1pt height5.3pt depth 0pt\relax
            \vrule width 1pt height5.7pt depth 0pt\relax
            \vrule width 1pt height6.0pt depth 0pt\relax
            \vrule width 1pt height6.3pt depth 0pt\relax
            \vrule width 1pt height6.7pt depth 0pt\relax
            \vrule width 1pt height7.0pt depth 0pt\relax
            \vrule width 1pt height7.3pt depth 0pt\relax
            \vrule width 1pt height7.7pt depth 0pt\relax
            \vrule width 1pt height8.0pt depth 0pt\relax
            \vrule width 1pt height8.3pt depth 0pt\relax
            \vrule width 1pt height8.7pt depth 0pt\relax
            \vrule width 1pt height9pt depth 0pt\relax
      }


\def\widestnumber#1#2{}

\def\citewarning#1{\ifx\shlhetal\relax 
    \else
    \par{#1}\par
    \fi
}

\def\rm{\fam0 \tenrm}

\def\fakesubhead#1\endsubhead{\bigskip\noindent{\bf#1}\par}



%
%
%

%

\font\textrsfs=rsfs10
\font\scriptrsfs=rsfs7
\font\scriptscriptrsfs=rsfs5

\newfam\rsfsfam
\textfont\rsfsfam=\textrsfs
\scriptfont\rsfsfam=\scriptrsfs
\scriptscriptfont\rsfsfam=\scriptscriptrsfs

\edef\oldcatcodeofat{\the\catcode`\@}
\catcode`\@11

\def\Cal@@#1{\noaccents@ \fam \rsfsfam #1}

\catcode`\@\oldcatcodeofat


\expandafter\ifx \csname margininit\endcsname \relax\else\margininit\fi

\long\def\red#1\endred{}
\long\def\green#1\endgreen{}
\long\def\blue#1\endblue{}

\def\endred{ \unmatched endred! }
\def\endgreen{ \unmatched endgreen! }
\def\endblue{ \unmatched endblue! }

\ifx\latexcolors\undefinedcs\def\latexcolors{}\fi

\def\emptycs{}
\def\evaluatelatexcolors{%
        \ifx\latexcolors\emptycs\else
        \expandafter\xxevaluate\latexcolors\xxfertig\evaluatelatexcolors\fi}
\def\xxevaluate#1,#2\xxfertig{\setupthiscolor{#1}%
        \def\latexcolors{#2}}

\font\smallfont=cmsl7
\def\rutgerscolor{\ifmmode\else\endgraf\fi\smallfont
\advance\leftskip0.5cm\relax}
\def\setupthiscolor#1{\edef\tmptmpcs{\noexpand\bgroup\noexpand\rutgerscolor
\noexpand\def\noexpand\currentcolor{#1}%
\noexpand}%
\expandafter\let\csname#1\endcsname\tmptmpcs
\def\tmptmpcs{\checkColorUnmatched{#1}\popthecolor}
\expandafter\let\csname end#1\endcsname\tmptmpcs}

\def\checkColorUnmatched#1{\def\expectcolor{#1}%
    \ifx\expectcolor\currentcolor   
    \else \edef\failhere{\noexpand\tryingToClose '\currentcolor' with end\expectcolor}\failhere\fi}

\def\currentcolor{???}

\def\popthecolor{\ifmmode\else\endgraf\fi\egroup}

\expandafter\def\csname#1\endcsname{}

\evaluatelatexcolors

 \let\outerhead\head
 \def\head{\innerhead}
 \let\innerhead\outerhead

 \let\outersubhead\subhead
 \def\subhead{\innersubhead}
 \let\innersubhead\outersubhead

 \let\outersubsubhead\subsubhead
 \def\subsubhead{\innersubsubhead}
 \let\innersubsubhead\outersubsubhead

 \def\proclaim{\innerproclaim}
 \let\innerproclaim\outerproclaim

 %
 %
 %
 %

\def\demo#1{\medskip\noindent{\it #1.\/}}
\def\enddemo{\smallskip}

\def\remark#1{\medskip\noindent{\it #1.\/}}
\def\endremark{\smallskip}

\pageheight{8.5truein}
\topmatter
\title{The depth of ultraproducts of Boolean Algebras} \endtitle
\author {Saharon Shelah \thanks {\null\newline I would like to thank 
Alice Leonhardt for the beautiful typing. \null\newline
This research was supported by the United States-Israel Binational
Science Foundation. Publication 853.} \endthanks} \endauthor 

\affil{The Hebrew University of Jerusalem \\
Einstein Institute of Mathematics \\
Edmond J. Safra Campus, Givat Ram \\
Jerusalem 91904, Israel
 \medskip
 Department of Mathematics \\
 Hill Center-Busch Campus \\
  Rutgers, The State University of New Jersey \\
 110 Frelinghuysen Road \\
 Piscataway, NJ 08854-8019 USA} \endaffil

\abstract   We show in ZFC, that the depth of ultraproducts of Boolean
Algebras may be bigger than the ultraproduct of the depth of those
Boolean Algebras. \endabstract
\endtopmatter
\document

\newpage

\head {\S0 Introduction} \endhead  \resetall \sectno=0
 \spuriousreset
\bigskip

Monk has looked systematically at cardinal invariants of Boolean
Algebras.  In particular, he has looked at the relations between
inv$(\dsize \prod_{i < \kappa} \bold B_i/D)$ and $\dsize \prod_{i <
\kappa}$ inv$(\bold B_i)/D$, i.e., the invariant of the ultraproducts
of a sequence of Boolean Algebras vis the ultraproducts of the sequence
of the invariants of those Boolean Algebras 
for various cardinal invariants inv of
Boolean Algebras.  That is: is it always true that inv$(\dsize \prod_{i
< \kappa} \bold B_i/D) \le \text{ inv}(\dsize \prod_{i < \kappa} \bold
B_i/D)$? is it consistently always true?  Is it always true that
$\dsize \prod_{i < \kappa}$ inv$(\bold B_i)/D \le 
\text{ inv}(\Pi \bold B_i/D)$? is it
consistenly always true?  See more on this in Monk \cite{Mo96}.
Roslanowski Shelah \cite{RoSh:534} deals with specific inv and with more
on kinds of cardinal invariants and
their relationship with ultraproducts.  Monk \cite{Mo90a}, \cite{Mo96},
in his list of open problems raises the question for the central
cardinal invariants, most of them have been solved by now; see Magidor
Shelah \cite{MgSh:433}, Peterson \cite{Pe97}, Shelah \cite{Sh:345},
\cite{Sh:462}, \cite{Sh:479}, \cite[\S4]{Sh:589}, \cite{Sh:620}, \cite{Sh:641},
\cite{Sh:703}, Shelah and Spinas \cite{ShSi:677}.

We here solve problem 12 of \cite{Mo96}, pg.287 in ZFC
constructing an example.  This example works for the length too.  As
in several earlier cases we use pcf theory to resolve the question
near singular cardinals, see \cite{Sh:g}.
\newpage

\head {\S1 On problem 12,p.287 of \cite{Mo96},Monk 1996 book} \endhead  \resetall \sectno=1
 \spuriousreset
\bigskip

\proclaim{\stag{1.1} Claim}  1) Assume
\mr
\item "{$(a)$}"  $\mu = \mu^\kappa > 2^\kappa$
\sn
\item "{$(b)$}"  $\mu$ singular, {\rm cf}$(\mu) = \theta$.
\ermn
\ub{Then} there are Boolean Algebras $B_i$ for $i < \kappa$
such that
\mr
\item "{$(\alpha)$}"  {\rm Depth}$(\bold B_i) \le \mu$ for $i < \kappa$,
hence
\sn
\item "{$(\alpha)'$}"  for any ultrafilter $D$ on $\kappa,\mu = 
\dsize \prod_{i < \kappa} { \text{\rm Depth\/}}(\bold B_i)/D$
\sn
\item "{$(\beta)$}"  for any uniform ultrafilter $D$ on $\kappa$, the
Boolean Algebra $\dsize \prod_{i < \kappa} B_i/D$ has depth $\ge \mu^+$.
\ermn
2) We can replace in $(\alpha) + (\beta)$, $\mu$ by $\mu_1$ if $\mu_1
< { \text{\rm pp\/}}(\mu)$ except in very rare cases, in particular it
suffices to assume that {\rm
pp}$_{J^{\text{bd}}_{\text{cf}(\mu)}}(\mu) 
> \mu_1$ \ub{or} it suffices
then $\lambda = { \text{\rm tcf\/}}(\prod a,<_J),J$ an ideal of ${\frak
a},\emptyset = \dbca_{i < \kappa} {\frak b}_i,{\frak b}_i \in J$
decreasing with empty intersections $\theta \in {\frak a} \Rightarrow
{ \text{\rm max pcf\/}}({\frak a} \cap \theta) < \lambda$.
\endproclaim
\bigskip

\demo{Proof of \scite{1.1}}  This is a special case of \scite{1.3}.
\enddemo
\bigskip

\remark{\stag{1.2} Remark}  Clearly for any given $\kappa$ there are
many such $\mu$'s, e.g., $\beth_{\kappa^+}$.
\endremark
\bigskip

\proclaim{\stag{1.3} Claim}  Assume
\mr
\item "{$(a)$}"  $J$ is an ideal on ${\frak a}$, {\rm sup}$({\frak a}) =
\mu,\mu$ singular and $\mu = \lim_J({\frak a})$; that is
$(\forall \mu_1 < \mu)[{\frak a} \cap \mu_1 \in J)$
and $\theta \in {\frak a} \Rightarrow$ {\rm max pcf}$({\frak a} \cap \beta) <
\mu$
\sn
\item "{$(b)$}"  $\lambda = { \text{\rm tcf\/}}(\prod {\frak a},\le_J)$
as witnessed by $\langle f_\alpha:\alpha < \lambda \rangle$
\sn
\item "{$(c)$}"  ${\frak b}_i \in J^+$ for $i < \kappa,{\frak b}_i$
decreasing with $i$ and $\emptyset = \cap\{{\frak b}_i:i < \kappa\}$
\sn  
\item "{$(d)$}"  $D$ is a uniform ultrafilter on $\kappa$.
\ermn
\ub{Then} for some sequence $\langle \bold B_i:i < \kappa \rangle$ of
Boolean Algebras we have:
\mr
\item "{$(\alpha)$}"  {\rm Depth}$^+(\dsize \prod_{i < \kappa} \bold B_i/D) >
\lambda$ (if $\lambda = \mu^+$ this means 
{\rm Depth}$(\dsize \prod_{i < \kappa} \bold B_i/D) \ge \lambda)$
\sn
\item "{$(\beta)$}"  {\rm Depth}$^+(\bold B_i) \le \lambda$ 
(if $\lambda = \mu^+$ this means {\rm Depth}$(\bold B_i/D) \le \mu)$.
\endroster
\endproclaim
\bigskip

\demo{Proof}  We can find $\langle f_\alpha:\alpha < \lambda \rangle$
which is $<_J$-increasing cofinal in $(\pi {\frak a},<_J)$ and
satisfies $\theta \in {\frak a} \Rightarrow |\{f_\alpha \restriction
({\frak a} \cap \theta):\alpha < \lambda\}| < \theta$ (see
\cite[II,3.5,pg.65]{Sh:g}).  
We define a function $\theta:[\lambda]^2 \rightarrow \kappa$ by: for 
$\alpha \ne \beta < \lambda$ we let $\theta\{\alpha,\beta\} 
= \text{ Min}\{\theta \in {\frak
a}:f_\alpha(\theta) \ne f_\beta(\theta)\}$ and we define a two place
relation $<_i$ on $\lambda$ by: $\alpha <_i \beta$ \ub{iff} 
$\theta\{\alpha,\beta\} \in {\frak a} \backslash {\frak b}_i \and
f_\alpha(\theta\{\alpha,\beta\}) < f_\beta(\theta\{\alpha,\beta\})$.
Now
\mr
\item "{$\circledast_1$}"   $\le_i$ is a partial order of $\lambda$.
\nl
[Why?  Assume $\alpha <_i \beta <_i \gamma$.
\ermn
Now
\mn
\ub{Case 1}: $\alpha = \beta \vee \beta = \gamma$: trivial.
\mn
\ub{Case 2}:   $\theta\{\alpha,\beta\} < \theta\{\beta,\gamma\}$ so

$$
\theta\{\alpha,\gamma\} = \theta\{\alpha,\beta\},
$$

$$
\align
(f_\alpha(\theta\{\alpha,\beta\}),f_\beta(\theta\{\alpha,\beta\}) &=
 (f_\alpha(\theta\{\alpha,\gamma\}),f_\beta(\theta\{\alpha,\gamma\})  \\
  &= (f_\alpha(\theta\{\alpha,\gamma\}),f_\gamma(\theta\{\alpha,\gamma\})
\endalign
$$
\mn
and we are done.
\mn
\ub{Case 3}:  $\theta\{\alpha,\beta\} > \theta\{\beta,\gamma\})$.

Similarly.
\mn
\ub{Case 4}:  $\theta(\alpha,\beta) = \theta(\beta,\gamma)$.  Call it
$\theta$.  So $f_\alpha \restriction \theta = f_\beta \restriction
 \theta = f_\gamma \restriction \theta$ and $f_\alpha(\theta) <
 f_\beta(\theta) < f_\gamma(\theta)$, hence $\theta\{\alpha,\gamma\} =
 \theta$ and $f_\alpha <_i f_\gamma$ as required.  So $\circledast_1$
 holds.]

Let $\bold B_i = BA[(\lambda,<_i)]$ where for a partial order
$(I,\le_I),BA[(I,\le_I)]$ is the Boolean Algebra generated by $\{x_t:t
< I\}$ freely except that
\mr
\item "{$\circledast_2$}"  $x_s \le x_t$ when $s \le_I x_t$.
\ermn
Now
\mr
\item "{$\circledast_3$}"  in $\bold B = \dsize \prod_{i < \kappa} 
\bold B_i/D$, there is an increasing sequence of length $\lambda$.
\nl
[Why?  Let $a_\alpha = \langle x_\alpha:\alpha < \kappa \rangle/D$,
now if $\alpha < \beta$ then $\theta\{\alpha,\beta\} \notin {\frak
b}_i \Rightarrow B_i \models ``x_\alpha < x_\beta"$ and $\alpha <
\beta \Rightarrow \{i < \kappa:\theta(\theta,\beta) \notin {\frak
b}_i\} \in D$ as $D$ is a uniform ultrafilter on $\kappa$ and the
sequence $\langle {\frak b}_i:i < \kappa \rangle$
decreases with intersection $\emptyset$ we are done easily.  Together
$\alpha < \beta \Rightarrow B \models a_\alpha < a_\beta$; so $\langle
a_\alpha;\alpha < \lambda \rangle$ is as required so $\circledast_3$ holds.]
\ermn
So it is enough to prove (as done in the rest of the proof).
\mr
\item "{$\circledast_4$}"  Depth$^+(\bold B_i) \le \lambda$.
\ermn
Toward contradiction, assume $\langle a_\alpha:\alpha < \lambda
\rangle$ is an $<_{B_i}$-increasing sequence of members of $\bold B_i$.  Let
$a_\alpha =
\sigma_\alpha(x_{\gamma(\alpha,0)},\dotsc,x_{(\alpha,n_\alpha-1)})$
where $\sigma_\alpha$ is a Boolean term and

$$
\gamma(\alpha,0) < \gamma(\alpha,1) < \ldots <
\gamma(\alpha,n_\alpha-1) < \lambda.
$$
\mn
Without loss of generality $\sigma_\alpha = \sigma_*$ so $n_\alpha =
n_*$ and $\theta\{\gamma(\alpha,\ell_1),\gamma(\alpha,\ell_2)\}$ is
the same for all $\alpha < \lambda$, say is $\theta_{\ell_1,\ell_2}$.
Now \wilog \, for some $\theta_* \in {\frak a}$ satisfying $\ell_1 < \ell_2 <
n_* \Rightarrow \theta_{\ell_1,\ell_2} <\theta_*$ 
we have $\ell < n_* \and \alpha < \beta <
\lambda \Rightarrow f_{\gamma(\alpha,\ell)} \restriction ({\frak a}
\cap \theta_*) = f_{\gamma(\beta,\ell)} \restriction ({\frak a} \cap
\theta_*)$. \nl
[Why?  Recall that $\theta \in {\frak a} \Rightarrow \theta >
|\{f_\alpha \restriction \theta:\alpha < \lambda\}|$ and $\theta \in
{\frak a} \Rightarrow \theta < \lambda = \text{ cf}(\lambda)$.]

Also \wilog \, for some $m_* < n_*,\ell < m_* \Rightarrow
\gamma(\alpha,\ell) = \gamma_*(\ell)$ and $\alpha < \beta \Rightarrow
\gamma(\alpha,n_*-1) < \gamma(\beta,m_*)$.  By
\cite[II,4.10A,4.10B,pg.76,77]{Sh:g} as 
${\frak b}_i \in J^+$ we can find $\theta^* \in {\frak b}_i \backslash
\theta_*$ and $\alpha < \beta$ such that
\mr
\item "{$\boxtimes^\theta_{\alpha,\beta}$}"  $\theta^* =
\theta\{\gamma(\alpha,\ell_1),\gamma(\beta,\ell_2)\}$ whenever
$\ell_1 \ne \ell_2 \in \{m_*,\dotsc,n_*-1\}$.
\ermn
Now let $I =
\{\gamma(\alpha,\ell),\gamma(\beta,\ell):\ell < n_*\}$.  Now we know that
$BA[(I,\le_i \restriction I)]$ is a Boolean subalgebra of $B_i$ hence
$BA[(I,\le_i \restriction I)] \models
``\sigma_*(\ldots,x_{\gamma(\alpha,\ell)},\ldots) <
\sigma_*(\ldots,x_{\gamma(\beta,\ell)},\ldots)"$.  But every
automorphism $\pi$ of $(I,\le_i \restriction I)$ induces an automorphism
$\hat \pi$ of $BA[(I,\le_i \restriction I)]$, but the permutation
$\pi$ interchanging $\gamma(\alpha,\ell)$ with $\gamma(\beta,\ell)$ is
an automorphism of $(I,\le_i \restriction I)$ so $BA[(I,\le_i
\restriction I)] \models ``\hat
\pi(\sigma_x(\ldots,x_{\gamma(\alpha,\ell)},
\ldots) < \hat \pi(\sigma_*(\ldots,x_{\gamma(\beta,\ell)},\ldots))"$ 
but this gives a contradiction.  \hfill$\square_{\scite{1.3}}$\margincite{1.3}
\enddemo
\bigskip

\proclaim{\stag{1.4} Claim}  In \scite{1.3} we can add
\mr
\item "{$(\alpha)''$}"  Length$^+(\bold B_i) \le \lambda$ for $i <
\kappa$.
\endroster
\endproclaim
\bigskip

\demo{Proof}  The same proof works, only concerning $\circledast_4$,
it is now
\mr
\item "{$\circledast_i$}"  Length$(\bold B_i) \le \lambda$.
\ermn
The proof is the same but we do not know that BA$[(I,\le_i
\restriction I)] \models
``\sigma_*(\ldots,x_{\gamma(\alpha,\ell)},\ldots)_{\ell < n_*} <
\sigma_*(\ldots,x_{\gamma(\beta,\ell)},\ldots)_{\ell < n_*}"$ but only
know that BA$[(I,\le_i \, \restriction I)] \models$ ``the elements
$\sigma_*(\ldots,x_{\gamma(\alpha,\ell)},\ldots)_{\ell < n_*}$ and
$\sigma_*(\ldots,x_{\gamma)\beta,\ell)},\ldots)_{\ell < n_*}$ are
comparable. \hfill$\square_{\scite{1.4}}$\margincite{1.4}
\enddemo
\bigskip

\nocite{ignore-this-bibtex-warning} 
\newpage
    
REFERENCES.  
\bibliographystyle{lit-plain}
\bibliography{lista,listb,listx,listf,liste}

\enddocument